\documentclass[smallcondensed]{svjour3}     
\usepackage{graphicx}
\usepackage[utf8]{inputenc}
\usepackage[english]{babel}
\usepackage{amsmath}
\usepackage{amsfonts}
\usepackage{textcomp}
\usepackage{multirow}
\usepackage{color}
\usepackage{subfigure}

\newtheorem{proposizione}{Proposition}

\newcommand{\ve}[1]{\boldsymbol{#1}}%

\begin{document}

\title{A new steplength selection for scaled gradient methods with application to image deblurring}

\titlerunning{Non-negative image deblurring: a limited memory approach}

\author{Federica Porta \and Marco Prato	\and Luca Zanni}

\institute{F. Porta \and M. Prato \and L. Zanni
					 \at
					 Dipartimento di Scienze Fisiche, Informatiche e Matematiche, Universit\`{a} degli Studi di Modena e Reggio Emilia, Via Campi 213/b, 41125 Modena, Italy\\
              Tel.: +39-059-2055590\\
              Fax: +39-059-2055216\\
              \email{marco.prato@unimore.it}}

\date{Received: date / Accepted: date}

\maketitle

\begin{abstract} 
Gradient methods are frequently used in large scale image deblurring problems since they avoid the onerous computation of the Hessian matrix of the objective function. Second order information is typically sought by a clever choice of the steplength parameter defining the descent direction, as in the case of the well-known Barzilai and Borwein rules. In a recent paper, a strategy for the steplength selection approximating the inverse of some eigenvalues of the Hessian matrix has been proposed for gradient methods applied to unconstrained minimization problems. In the quadratic case, this approach is based on a Lanczos process applied every $m$ iterations to the matrix of the gradients computed in the previous $m$ iterations, but the idea can be extended to a general objective function. In this paper we extend this rule to the case of scaled gradient projection methods applied to constrained minimization problems, and we test the effectiveness of the proposed strategy in image deblurring problems in both the presence and the absence of an explicit edge-preserving regularization term.
\end{abstract}

\section{Problem formulation}

The image formation process is an inverse problem that can be modeled as the following linear system
\begin{equation}
\label{Linear_System}
\ve{y} = A\ve{x} + \ve{b} + \ve{\eta}\,,
\end{equation}
where $\ve{y} \in \mathbb{R}^{n^2}$ is the non-negative observed data, $\ve{x}\in \mathbb{R}^{n^2}$ represents an ideal, undistorted image to be recovered, $A \in \mathbb{R}^{n^2 \times n^2}$ is a typically ill-conditioned matrix describing the blurring effect, $\ve{b} \in \mathbb{R}^{n^2}$ is a known non-negative background radiation and $\ve{\eta} \in \mathbb{R}^{n^2}$ is the noise corrupting the data. A typical assumption for the matrix $A$ is that it has non-negative elements and each row and column has at least one positive entry. Because of the ill-conditioning affecting the problem and the presence of noise on the measured data, a trivial approach that seeks the solution of \eqref{Linear_System} is in general not successful; thus, alternative strategies must be exploited. Variational approaches to image restoration \cite{Bertero2008,Vogel2002} suggest to recover the unknown object through iterative schemes suited for the following constrained minimization problem
\begin{equation}
\label{Min_prob}
\min_{\ve{x} \geq \ve{0}} J_0(\ve{x})
\end{equation}
where $J_0$ is a continuously differentiable convex function measuring the difference between the model and the data. The definition of the function $J_0$ depends on the noise type introduced by the acquisition system. Particularly, in the case of additive white Gaussian noise the cost function is characterized by a least squares distance of the form
\begin{equation}
\label{J_Gaussian}
J_0(\ve{x}) = J_0^{LS}(\ve{x}) = \frac{1}{2} \| A\ve{x} + \ve{b} - \ve{y}\|^2 \, ,
\end{equation}
while, when the data are affected by Poisson noise, the so-called Kullback-Leibler (KL) divergence is used:
\begin{equation}
\label{J_Poisson}
J_0(\ve{x}) = J_0^{KL}(\ve{x}) = \sum_{i=1}^{n^2} \left\{ y_i \ln\frac{y_i}{(A\ve{x} + \ve{b})_i} + (A\ve{x} + \ve{b})_i - y_i\right\} \, ,
\end{equation}
where we assume that $0\ln 0 = 0$ and $(A\ve{x} + \ve{b})_i > 0$, $\forall i=1,\ldots,n^2$. In both cases, taking into account also the assumptions on $A$, we may observe that $J_0$ is non-negative, convex and coercive on the non-negative orthant, which means that problem \eqref{Min_prob} has global solutions. Moreover, if the equation $A\ve{x}=\ve{0}$ has only the solution $\ve{x}=\ve{0}$, then $J_0^{LS}$ is strictly convex, while the same conclusion holds for $J_0^{KL}$ if the additional condition $y_i > 0$, $\forall i=1,\ldots,n^2$, is satisfied \cite{Bertero2010,Bonettini2009}. In these settings, the strict convexity of $J_0$ implies that the solution of \eqref{Min_prob} is unique.\\
Due to the ill-posedness of the image restoration problem, one is not interested in computing the minimum points of $J_0$ in \eqref{J_Gaussian} or \eqref{J_Poisson} because the exact solution of \eqref{Min_prob} does not provide a sensible estimate of the unknown image. For this reason, iterative minimization methods are usually exploited to obtain acceptable solutions by arresting the algorithm before convergence through some stopping criteria, as the classic Morozov's discrepancy principle in the case of Gaussian noise (see e.g. \cite{Hansen1997}) or some recently proposed strategies for Poisson data \cite{Bardsley2009,Bertero2010,Carvalan2011}.\\
Another technique to tackle to this problem requires to exactly solve the following optimization problem
\begin{equation}
\label{Min_prob_reg}
\min_{\ve{x} \geq \ve{0}} J_0(\ve{x}) + \beta J_R(\ve{x}),
\end{equation}
where $J_R$ is a regularization term adding a priori information on the solution and $\beta$ is a positive parameter balancing the role of the two objective function components $J_0$ and $J_R$. A frequently used function for the regularization term is a smooth approximation of the total variation, also known in the literature as {\em hypersurface potential} (HS), defined as \cite{Acar2004,Bertero2010}
\begin{equation}
\label{TV_smooth}
 J_R(\ve{x}) =
J_R^{HS}(\ve{x}) = \sum_{i,j=1}^n \sqrt{((\mathcal{D} \ve{x})_{i,j})_1^2+((\mathcal{D} \ve{x})_{i,j})_2^2+\delta^2},
\end{equation}
where the discrete gradient operator $\mathcal{D}: \mathbb{R}^{n^2}\longrightarrow\mathbb{R}^{2n^2}$ is set through the standard finite difference with periodic boundary conditions
\begin{equation}\label{discr_grad}
(\mathcal{D} \ve{x})_{i,j} =
\left(
\begin{array}{c}
 ((\mathcal{D}\ve{x})_{i,j})_1\\
 ((\mathcal{D}\ve{x})_{i,j})_2
\end{array}
\right)=
\left(
\begin{array}{c}
 x_{i+1,j}-x_{i,j}\\
 x_{i,j+1}-x_{i,j}\\
\end{array}
\right), \ \ \ x_{n+1,j} = x_{1,j}, \ \ \ x_{i,n+1} = x_{i,1}.
\end{equation}
When $J_R=J_R^{HS}$ and $J_0$ is one of the two considered cost functions, the objective function in \eqref{Min_prob_reg} is non-negative, strictly convex and coercive on the non-negative orthant \cite{Bertero2010}. It follows that problem \eqref{Min_prob_reg} has a unique solution.\\
Both formulations of the imaging problem require an effective optimization method able to provide a meaningful solution in a reasonable time. Among all possible choices, first-order methods are particularly suited to deal with this kind of problems for several reasons. First, due to the large size of the images (which becomes a crucial issue especially in 3D applications), the handling of the Hessian matrix is an impractical task. Then, first-order methods are used to quickly achieve solutions with low/medium accuracy, which is a general requirement in imaging problems. Finally, when the optimization scheme is used as iterative regularization method to minimize the cost function \eqref{Min_prob}, an excessively fast convergence makes the automatic choice of the stopping iteration a crucial issue, since a difference of few iterations from the one providing the best reconstruction can lead to substantial differences in the final images.\\
In this paper we extend to the case of a general scaled gradient projection method \cite{Bertsekas1999,Birgin2003,Bonettini2009} a steplength selection rule recently proposed by Fletcher \cite{Fletcher2012} in the unconstrained optimization framework and we test its effectiveness in image deblurring problems. This rule is based on the estimate of some eigenvalues of the Hessian matrix which, for quadratic problems, can be achieved by means of a Lanczos process applied to a certain number of consecutive gradients. Since the scheme depends only on these stored gradients, it can be generalized to nonquadratic objective functions, showing very competitive results in several benchmark problems with respect to other first-order and quasi-Newton methods. The extension to scaled gradient projection methods applied to non-negatively constrained problems requires a generalization of the matrix with the last gradients accounting for the presence of both the scaling matrix multiplying the gradient and the projection on the non-negative orthant. The resulting scheme consists in the storage of a set of scaled gradients (instead of the usual ones) in which some components of the gradients themselves are put equal to zero. Our numerical experiments on the non-negative minimization of the LS distance and the KL divergence show that the proposed approach is able to compete with standard gradient methods and other recently proposed schemes, providing in some cases good reconstructions with a significantly lower number of iterations.\\
The plan of the paper is the following: in section \ref{sec2} we recall the features of a scaled gradient projection method and, in particular, of the scaling matrix multiplying the gradient. In section \ref{sec3} we focus the analysis on the choice of the steplength parameter and we describe state-of-the-art strategies and our proposed rule. In section \ref{sec4} some numerical experiments on small quadratic programming (QP) and image deblurring least-squares problems are presented, while in section \ref{sec5} we address the image deblurring problem with data perturbed with Poisson noise also by adding an edge-preserving regularization term in the objective function. Some ideas on a possible generalization of the proposed rule to different constraints are provided in section \ref{sec6}, together with a numerical test on the Rudin-Osher-Fatemi model \cite{Rudin1992}. Our conclusions are given in section \ref{sec7}.

\section{Scaled gradient projection methods}\label{sec2}

A general scaled gradient projection (SGP) method \cite{Bonettini2009} for the solution of
\begin{equation}
\label{min_prob_gen}
\min_{\ve{x}\geq 0} J(\ve{x}),
\end{equation}
with $J$ differentiable function, is an iterative algorithm whose $(k+1)$-th iteration is defined by
\begin{equation}
\label{SGP_iter}
\ve{x}^{(k+1)} = \ve{x}^{(k)} + \lambda_k \ve{d}^{(k)} = \ve{x}^{(k)} + \lambda_k\left(\mathbb{P}_{+,D_k^{-1}}(\ve{x}^{(k)} - \alpha_k D_k\ve{g}^{(k)}) -\ve{x}^{(k)}\right),
\end{equation}
where
\begin{itemize}
\item[$\bullet$]  $\ve{x}^{(k)}\ge \ve{0}$;
\item[$\bullet$] $\ve{g}^{(k)} = \nabla J(\ve{x}^{(k)})$ is the gradient of the objective function at iteration $\ve{x}^{(k)}$;
\item[$\bullet$] $\lambda_k \in (0,1]$ is a linesearch parameter ensuring a sufficient decrease of the objective function along the descent direction $\ve{d}^{(k)}$, e.g. by means of an Armijo rule \cite{Bertsekas1999};
\item[$\bullet$] $\alpha_k$ is a positive steplength chosen in a fixed range $[\alpha_{\min},\alpha_{\max}]$, with $0 < \alpha_{\min} < \alpha_{\max}$;
\item[$\bullet$] $D_k$ is a symmetric and positive definite scaling matrix with eigenvalues lying in a fixed positive interval $[L_1,L_2]$;
\item[$\bullet$] $\mathbb{P}_{+,D}(\cdot)$ denotes the projection operator onto the non-negative orthant with respect to the norm induced by the matrix $D$:
\begin{equation*}
\mathbb{P}_{+,D}(\ve{x}) = {\rm{arg}}\min_{\ve{y}\geq 0} \|\ve{y}-\ve{x}\|_D = {\rm{arg}}\min_{\ve{y}\geq 0} \frac{1}{2}\ve{y}^TD\ve{y}-\ve{y}^TD\ve{x}.
\end{equation*}
\end{itemize}
The boundedness conditions on the steplengths and the eigenvalues of the scaling matrices are necessary to prove the convergence result for this method (see \cite[Theorem 2.1]{Bonettini2009}), that we report for completeness.
\begin{theorem}
 Let $\ve{x}^{(0)}\ge \ve{0}$ and
assume that the level set $\Omega_0 = \{\ve{x}\geq\ve{0} : J(\ve{x}) \leq J(\ve{x}^{(0)})\}$ is bounded. Every accumulation point of the sequence $\{\ve{x}^{(k)}\}_{k\in\mathbb{N}}$ generated by the SGP algorithm is a stationary point of \eqref{min_prob_gen}.
\end{theorem}
When SGP is applied to the imaging minimization problem \eqref{Min_prob} or \eqref{Min_prob_reg}, the coercivity of the objective function on the non-negative orthant assures that $\Omega_0$ is bounded for any $\ve{x}^{(0)} \geq \ve{0}$, therefore the sequence generated by SGP is bounded and admits limit points;
the uniqueness of the limit point is ensured when the objective function is strictly convex.\\
In imaging applications, the scaling matrix $D_k$ is usually chosen according to the cost function $J_0$ and the regularization term $J_R$. Following the approach proposed in \cite{Lanteri2002,Lanteri2001}, if $\nabla J_0$ and $\nabla J_R$ can be decomposed in the form
\begin{equation}
\label{grad_dec}
-\nabla J_0(\ve{x}) = U_0(\ve{x}) - V_0(\ve{x}) \quad ; \qquad -\nabla J_R(\ve{x}) = U_R(\ve{x}) - V_R(\ve{x}),
\end{equation}
with $U_0,U_R \geq 0$ and $V_0,V_R > 0$, then a possible scaling matrix is given by
\begin{equation}\label{Dk}
D_k = \max\left( L_1, \min\left( L_2, {\rm{diag}}\left(\frac{\ve{x}^{(k)}}{V_0(\ve{x}^{(k)}) + \beta V_R(\ve{x}^{(k)})}\right)\right)\right), \qquad L_1 \le L_2,
\end{equation}
where $\ve{v}/\ve{w}$ is the componentwise ratio between $\ve{v}$ and $\ve{w}$. We remark that the choice of a diagonal scaling matrix is preferable since in this case the projection on the non-negative orthant is straightforward and does not require the solution of a further quadratic subproblem at each iteration. Since in general the imaging matrix $A$ has non-negative entries, the gradients of the cost functions in \eqref{J_Gaussian} and \eqref{J_Poisson} satisfy the decomposition in \eqref{grad_dec}
\begin{equation*}
-\nabla J_0^{LS}(\ve{x}) = \underbrace{A^T\ve{y}}_{U_0^{LS}(\ve{x})} - \underbrace{A^T(A\ve{x} + \ve{b})}_{V_0^{LS}(\ve{x})} \quad ; \quad -\nabla J_0^{KL}(\ve{x}) = \underbrace{A^T\frac{\ve{y}}{A\ve{x} + \ve{b}}}_{U_0^{KL}(\ve{x})} - \underbrace{A^T\ve{1}}_{V_0^{KL}(\ve{x})},
\end{equation*}
where $\ve{1}$ is the vector with all entries equal to 1. In a similar way, the negative gradient of the regularization term in \eqref{TV_smooth} can be written as in \eqref{grad_dec} with \cite{Zanella2009}
\begin{align*}
[U_R^{HS}(\ve{x})]_{i,j} \! =&  \frac{x_{i+1,j} + x_{i,j+1}}{\sqrt{((\mathcal{D}\ve{x})_{i,j})_1^2 \! + \! ((\mathcal{D}\ve{x})_{i,j})_2^2+\delta^2}} \!
                              + \! \frac{x_{i,j-1}}{\sqrt{((\mathcal{D}\ve{x})_{i,j-1})_1^2 \! + \! ((\mathcal{D}\ve{x})_{i,j-1})_2^2 \! + \! \delta^2}} \\
												     &+ \frac{x_{i-1,j}}{\sqrt{((\mathcal{D}\ve{x})_{i-1,j})_1^2+((\mathcal{D}\ve{x})_{i-1,j})_2^2+\delta^2}},\\
[V_R^{HS}(\ve{x})]_{i,j} \! =&  \frac{2x_{i,j}}{\sqrt{((\mathcal{D}\ve{x})_{i,j})_1^2 \! + \! ((\mathcal{D}\ve{x})_{i,j})_2^2 \! + \! \delta^2}} \!
                              + \! \frac{x_{i,j}}{\sqrt{((\mathcal{D}\ve{x})_{i,j-1})_1^2 \! + \! ((\mathcal{D}\ve{x})_{i,j-1})_2^2 \! + \! \delta^2}} \\
												     &+ \frac{x_{i,j}}{\sqrt{((\mathcal{D}\ve{x})_{i-1,j})_1^2+((\mathcal{D}\ve{x})_{i-1,j})_2^2+\delta^2}}.
\end{align*}
The crucial task of speeding up the convergence of a scaled gradient projection method is generally assigned to the steplength parameter, which will be analyzed in the following section.

\section{A new steplength selection rule}\label{sec3}

Once the scaling matrix has been fixed, the steplength parameter $\alpha_k$ is chosen to encode some second order information to improve the converge rate of the scheme. Possible choices are the two rules proposed by Barzilai and Borwein (BB) \cite{Barzilai1988} for nonscaled gradient methods and extended by Bonettini et al \cite{Bonettini2009} to account for the presence of a scaling matrix $D_k$. These rules arise from the approximation of the Hessian $\nabla^2 J(\ve{x}^{(k)})$ with the diagonal matrix $B(\alpha_k)=(\alpha_kD_k)^{-1}$ and by imposing the following quasi-Newton properties on $B(\alpha_k)$:
\begin{align*}
\alpha_k^{BB1} &= \underset{\alpha_k \in \mathbb{R}}{\rm{argmin}} \|B(\alpha_k)\ve{s}^{(k-1)} - \ve{z}^{(k-1)}\|; \\
\alpha_k^{BB2} &= \underset{\alpha_k \in \mathbb{R}}{\rm{argmin}} \|\ve{s}^{(k-1)} - B(\alpha_k)^{-1}\ve{z}^{(k-1)}\|,
\end{align*}
where $\ve{s}^{(k-1)} = \ve{x}^{(k)} - \ve{x}^{(k-1)}$ and $\ve{z}^{(k-1)} = \nabla J(\ve{x}^{(k)}) - \nabla J(\ve{x}^{(k-1)})$. The resulting values become
\begin{equation*}
\alpha_k^{BB1} = \frac{{\ve{s}^{(k-1)}}^TD_k^{-2}\ve{s}^{(k-1)}}{{\ve{s}^{(k-1)}}^TD_{k}^{-1}\ve{z}^{(k-1)}} \qquad ; \qquad
\alpha_k^{BB2} = \frac{{\ve{s}^{(k-1)}}^TD_k\ve{z}^{(k-1)}}{{\ve{z}^{(k-1)}}^TD_{k}^{2}\ve{z}^{(k-1)}} \, ,
\end{equation*}
which reduce to the standard BB rules when $D_k$ is equal to the identity matrix $I$ for all $k$ (in the following, we will denote by GP a nonscaled gradient projection method). Many other steplength rules have been investigated in the last years (see \cite{Dai2003,DeAsmundis2013,DeAsmundis2014,Ruggiero2000,Yuan2006,Zhou2006} and references therein) and interesting convergence rate improvements have been observed by exploiting alternating criteria of the two BB rules, as the adaptive Barzilai-Borwein (ABB) method \cite{Zhou2006} and its generalizations ABB$_{\rm{min1}}$ and ABB$_{\rm{min2}}$ provided by Frassoldati et al \cite{Frassoldati2008}.\\
The aim of this paper is to realize an accelerating strategy for the SGP method through the generalization of a steplength selection rule recently suggested by Fletcher \cite{Fletcher2012} in the unconstrained optimization framework. For unconstrained minimization problems, theoretical considerations, confirmed by numerical experiments, showed the efficacy of this rule in improving the performances of first-order algorithms exploiting a single BB steplength rule. This analysis encouraged us to investigate the possibility of extending the Fletcher's scheme to the case of constrained optimization in order to use this innovative idea for scaled gradient projection method of the type \eqref{SGP_iter}, particularly in image deblurring applications. The new approach proposed in \cite{Fletcher2012} consists of a limited memory scheme defining the steplengths as the inverse of special approximations of the eigenvalue of the Hessian $\nabla^2 J(\ve{x})$. Let us consider a quadratic objective function $J(\ve{x}) = \frac{1}{2}\ve{x}^TA\ve{x}$, where $A$ is a symmetric and positive definite matrix. Then the steepest descent method applied to the unconstrained quadratic programming problem
\begin{equation}\label{unconst_QP}
\min_{\ve{x}\in \mathbb{R}^n} J(\ve{x})
\end{equation}
assumes the form
\begin{equation*}
\ve{x}^{(k+1)} = \ve{x}^{(k)} -\alpha_k\ve{g}^{(k)} =  \ve{x}^{(k)} -\alpha_k A\ve{x}^{(k)},\qquad k=0,1,\dots.
\end{equation*}
In particular, the following relation between the gradients holds true:
\begin{equation}\label{relgrad}
\ve{g}^{(k+1)} = \ve{g}^{(k)} -\alpha_k A\ve{g}^{(k)}.
\end{equation}
If a limited number $m$ of back values of the gradient vectors
\begin{equation}
\label{matrixG}
G = \left[\ve{g}^{(k-m)} \ \ldots \ \ve{g}^{(k-2)} \ \ve{g}^{(k-1)}\right]
\end{equation}
is stored in memory and the $(m+1)\times m$ matrix $\Gamma$ containing the reciprocals of the corresponding last $m$ steplengths is considered,
\begin{equation}
\Gamma=
\left[
\begin{array}{cccc}
\alpha^{-1}_{k-m}			&	&	&\\
-\alpha^{-1}_{k-m}	&\ddots		&	&\\
	&\ddots	&\alpha^{-1}_{k-2}		&\\
	&	&-\alpha^{-1}_{k-2}	&\alpha^{-1}_{k-1}\\
	&	&	&-\alpha^{-1}_{k-1}\\
\end{array}
\right],
\end{equation}
then equations \eqref{relgrad} for $k-m,\ldots,k-1$ can be rearranged in the matrix form
\begin{equation}\label{AG}
AG = [G \quad  \ve{g}^{(k)}]\Gamma.
\end{equation}
This equality can be used to rewrite the tridiagonal $m \times m$ matrix $\Phi$ provided by $m$ steps of the Lanczos iterative process applied to the matrix $A$ with starting vector $\ve{q}_1 = \ve{g}^{(k-m)}/\|\ve{g}^{(k-m)}\|$ \cite{Golub1996}. In fact, given an integer $m\ge 1$, the Lanczos process generates orthonormal vectors $\{\ve{q}_1, \ve{q}_2, \dots, \ve{q}_m\}$ that define a basis for the Krylov sequence $\{\ve{g}^{(k-m)}, A\ve{g}^{(k-m)}, \dots, A^{m-1}\ve{g}^{(k-m)} \}$ and such that the matrix
\begin{equation*}
\Phi = Q^TAQ,
\end{equation*}
where $Q=[\ve{q}_1, \ve{q}_2, \dots, \ve{q}_m]$, $Q^TQ = I$, is tridiagonal. Taking into account equation \eqref{relgrad} and that the columns of $G$ are in the space generated by the above Krylov sequence, we have $G=QR$, where $R$ is $m\times m$ upper triangular and nonsingular, assuming $G$ is full-rank. It follows from \eqref{AG} that the tridiagonal matrix $\Phi$ can be written as
\begin{equation*}
\Phi = Q^TAGR^{-1} = [R \quad Q^T\!\!\ve{g}^{(k)}]\Gamma R^{-1}
\end{equation*}
and, by introducing the vector $\ve{r} = Q^T \ve{g}^{(k)}$, that is the vector that solves the linear system $R^T \ve{r} = G^T\ve{g}^{(k)}$, we obtain
\begin{equation}\label{Phi}
\Phi = [R \quad \ve{r}]\Gamma R^{-1}.
\end{equation}
The eigenvalues of the tridiagonal matrix $\Phi$, called Ritz values, are approximations of $m$ eigenvalues of $A$ \cite{Golub1996} and, since $A$ is the Hessian matrix of the objective function $J$, they give some second order information about problem \eqref{unconst_QP}. The steplength selection rule proposed by Fletcher consists in exploiting the reciprocal of the $m$ Ritz values as steplengths in the next $m$ iterations. We refer to \cite{Fletcher2012} for a detailed motivation of this steplength rule and we focus on the features crucial for the extension of the rule to nonquadratic objective functions and to constrained optimization problems.
First of all we remark that \eqref{Phi} allows one to obtain the matrix $\Phi$ by simply exploiting the partially extended Cholesky factorization
\begin{equation*}
G^T [G \quad \ve{g}^{(k)}] = R^T [R \quad \ve{r}],
\end{equation*}
without the explicit use of the matrices $Q$ and $A$.
This is important both for the computational point of view and for the extension to nonquadratic functions.
For a general objective function, $\Phi$ is upper Hessenberg and the Ritz-like values are obtained by computing the eigenvalues of a symmetric and tridiagonal approximation $\widetilde{\Phi}$ of $\Phi$ defined as
\begin{equation*}
\widetilde{\Phi} = {\rm{diag}}(\Phi) + {\rm{tril}}(\Phi,-1) + {\rm{tril}}(\Phi,-1)^T,
\end{equation*}
where ${\rm{diag}}(\cdot)$ and ${\rm{tril}}(\cdot,-1)$ denote the diagonal and the strictly lower triangular parts of a matrix. Possible negative eigenvalues of the resulting matrix are discarded before using this set of steplengths for the next iterations. Several numerical experiments \cite{Fletcher2012}, for both quadratic and nonquadratic test problems, demonstrate that this new steplength selection rule is able to improve the convergence rate of steepest descent methods with respect to other, often used, possibilities for choosing the steplength.\\
Motivated by these promising results and taking into account that the convergence for the scaled gradient projection method \eqref{SGP_iter} is guaranteed for every choice of the steplength in a bounded interval, we tried to exploit the Fletcher's steplength selection rule in the algorithms used for constrained optimization. In the extension of the original scheme to the SGP method, the main change is the definition of a new matrix $\widetilde{G}$ that generalizes the matrix $G$ in \eqref{matrixG}. In particular, we have to consider two fundamental elements: the presence of the scaling matrix multiplying the gradient direction and the projection onto the feasible set. As concerns the former issue, we exploit the remark that each scaled gradient iteration can be viewed as a usual gradient iteration applied to a scaled objective function by means of a transformation of variables of the type $\ve{y} = D_{k}^{-1/2}\ve{x}$ \cite{Bertsekas1999}, where the notation $\displaystyle D^{1/2}$ indicates the square root matrix of $D$. This idea led us to store at each iteration the scaled gradient $D_{k}^{1/2}\ve{g}^{(k)}$ instead of $\ve{g}^{(k)}$. The non-negativity constraint is addressed by looking at the complementarity condition of the KKT optimality criteria \cite{Nocedal2006}, for which the components of the gradient related to inactive constraints in the solution have to vanish. To this aim, we emphasized the minimization over these components by storing the vectors $\widetilde{\ve{g}}^{(k)}$ whose $j$-th entry is given by
\begin{equation}\label{critnonneg}
\widetilde{g}^{(k)}_j = \begin{cases} 0 & {\rm{if}} \ x^{(k)}_j=0, \\ \left[\nabla J(\ve{x}^{(k)})\right]_j & {\rm{if}} \ x^{(k)}_j > 0. \end{cases}
\end{equation}
Driven by the previous considerations, our implementation of Fletcher's rule for the constrained case is based on the following choice for the matrix $\widetilde{G}$:
\begin{equation*}
\widetilde{G} = \left[D_{k-m}^{1/2}\widetilde{\ve{g}}^{(k-m)}, \ldots, D_{k-1}^{1/2}\widetilde{\ve{g}}^{(k-1)}\right].
\end{equation*}
As concerns the computational cost of the steplength derivation,  each group of $m$ iterations (called {\it sweep} in \cite{Fletcher2012}) requires the computation of the $m$ scaled gradients $D_{j}^{1/2}\widetilde{\ve{g}}^{(j)}$ and the $m \times m$ symmetric matrix $\widetilde{G}^T\widetilde{G}$, which can be performed with $m + (m+1)m/2 = (m+3)m/2$ vector-vector products. Since $m$ is typically a very small number (between 3 and 5), the Cholesky factorization of $\widetilde{G}^T\widetilde{G}$ and the solution of the linear system $R^T\ve{r}=\widetilde{G}^TD_k^{1/2}\widetilde{\ve{g}}^{(k)}$ are straightforward. It is worth noting that the computation of either the BB1 or the BB2 steplength for $m$ iterations needs 3$m$ vector-vector products. Therefore, if we assume for example $m=3$, then both the generalization of the limited memory approach and each BB steplength can be computed in $\mathcal{O}(9n^2)$ products, while the computational cost grows up to $\mathcal{O}(18n^2)$ for any alternating strategy of the two BB rules.\\
In the next sections we present the benefits that can be gained by using the steplength selection rule based on the Ritz values adapted to the constrained optimization in the image reconstruction framework.

\section{Numerical experiments - quadratic case}\label{sec4}

In this section we report the results of several numerical experiments we carried out on constrained QP problems in order to validate the efficacy of the limited memory selection rule. First we show few tests on the minimization of a quadratic function of 20 variables, with the analysis of the behaviour of three steplengths when varying some features of the optimization problem. Then we present realistic experiments of imaging problems with a comparison of several scaled and nonscaled gradient projection methods. All the numerical experiments have been performed by means of routines implemented by ourselves in Matlab$^\circledR$ R2010a and run on a PC equipped with a 1.60 GHz Intel Core i7 in a Windows 7 environment.

\subsection{Quadratic problems}\label{secG1}

The aim of this section is to investigate possible dependencies of the results provided by a (S)GP method with different steplengths on the features of the quadratic problem to be addressed, as the distribution of the eigenvalues of the Hessian matrix $A$, the number of active constraints and the condition number. Therefore, we built up some ad hoc tests to evaluate different selection rules for different choices of these parameters of the problem. In particular, we consider the minimization problem
\begin{equation}\label{toyQP}
\min_{\ve{x}\geq 0} \ve{x}^T A \ve{x} - \ve{y}^T \ve{x},
\end{equation}
where:
\begin{itemize}
\item[$\bullet$] we chose a vector $\ve{\xi}\in\mathbb{R}^{20}$ and we defined the matrix $A$ as $Q {\rm{diag}}(\ve{\xi}) Q^T$, where $Q$ is an orthogonal matrix obtained by a QR factorization of a random matrix;
\item[$\bullet$] we defined randomly the set $I_a \subseteq \{1,\ldots,20\}$ of $n_a$ active constraints;
\item[$\bullet$] we defined the vector of Lagrange multipliers $\ve{\mu}\in\mathbb{R}^{20}$ by setting $\mu_i=1$ if $i \in I_a$ and $\mu_i=0$ if $i \notin I_a$. In a similar way, we defined the solution of the problem $\ve{x}^*\in\mathbb{R}^{20}$ by setting $x^*_i=0$ if $i \in I_a$ and $x^*_i$ random in $(0,1)$ if $i \notin I_a$;
\item[$\bullet$] we defined the vector $\ve{y} = A\ve{x}^* - \ve{\mu}$.
\end{itemize}
The generalization of the limited memory (Ritz) steplength to the constrained case has been compared to the ABB$_{\rm{min1}}$ and BB1 values, where in the former case we used the generalized adaptive alternation rule
 proposed in \cite{Bonettini2009}. For all the three algorithms we exploited both a monotone and a nonmonotone linesearch \cite{Grippo1986} to determine the parameter $\lambda_k$. In the latter case, the sufficient decrease at each iteration is evaluated with respect to the maximum of the objective function on the last $M=10$ iterations. In the limited memory rule, the number $m$ of back stored gradient has been set equal to 3. Following \cite{Fletcher2012}, we started by considering $\ve{\xi}=((\sqrt{2})^0,\ldots,(\sqrt{2})^{19})$ and we investigated possible choices of the scaling matrix for the minimization problem \eqref{toyQP}. The number of active constraints has been set equal to 8. We remark that, since $A$ in our tests has also negative entries, the scaling matrix provided by the splitting of $\nabla J_0^{LS}$ in section 2 is not applicable. Possible scaling matrices are given by:
\begin{itemize}
\item[S1)] the inverse of the diagonal of $A$: $D_k^{PR} = {\rm{diag}}\left(\ve{1}/{\rm{diag}}(A)\right)$, which for the quadratic case is equivalent to apply a nonscaled gradient projection method to a preconditioned version of the minimization problem;
\item[S2)] the scaling matrix proposed by Coleman and Li \cite{Coleman1996} for interior trust region approaches applied to nonlinear minimization problems subject to box constraints: $D_k^{CL} = {\rm{diag}}\left(\widetilde{\ve{x}}^{(k)}\right)$, where $\widetilde{x}^{(k)}_i = x^{(k)}_i$ if $g^{(k)}_i \geq 0$ and $\widetilde{x}^{(k)}_i = 1$ if $g^{(k)}_i < 0$;
\item[S3)] the current iteration: $D_k^{XK} = {\rm{diag}}\left(\ve{x}^{(k)}\right)$.
\end{itemize}
The diagonal entries of all the scaling matrices have been projected in the range $[10^{-5},10^5]$ to guarantee the convergence of the schemes. In order to avoid the dependency of the analysis on the stopping criterion used, in Table \ref{tabG1} we reported the number of iterations required by the different algorithms to reach a relative reconstruction error (RRE) $\frac{\| \ve{x}^{(k)} -\ve{x}^*\|}{\|\ve{x}^*\|}$ lower than prefixed thresholds (e.g., $10^{-4}$, $10^{-6}$, $10^{-8}$). The performances with 
the trivial scaling matrix $D_k=I$ are also reported.\\

\begin{table}[ht]
\caption{Numbers of iterations required by SGP equipped with the limited memory (Ritz), ABB$_{\rm{min1}}$ and BB1 steplengths to reach RREs lower than $10^{-4}$, $10^{-6}$ and $10^{-8}$ for different scaling matrices (see text). The results obtained with a monotone ($M=1$) and nonmonotone ($M=10$) linesearch are reported. The asterisk denotes the maximum number of iterations allowed.}
\label{tabG1}
\begin{center}
\begin{tabular}{c|c|cc|cc|cc|}
\multicolumn{2}{c|}{}        & \multicolumn{2}{c|}{Ritz} & \multicolumn{2}{c|}{ABB$_{\rm{min1}}$} & \multicolumn{2}{c|}{BB1}\\
$D_k$                        & Tol       & M = 1  & M = 10 & M = 1 & M = 10 & M = 1  & M = 10 \\
\hline
\multirow{3}{*}{$I$}    		 & $10^{-4}$ &  100   & 94     &  136  & 93     &   155  & 124   \\
														 & $10^{-6}$ &  121   & 121    &  176  & 122    &   185  & 157   \\
														 & $10^{-8}$ &  161   & 127    &  219  & 155    &   230  & 178   \\
\hline
\multirow{3}{*}{$D_k^{PR}$}  & $10^{-4}$ &  90    & 82     &  120  & 114    &   143  & 124   \\
														 & $10^{-6}$ &  108   & 103    &  168  & 153    &   145  & 155   \\
														 & $10^{-8}$ &  157   & 115    &  222  & 163    &   220  & 193   \\
\hline
\multirow{3}{*}{$D_k^{CL}$}  & $10^{-4}$ &  282   & 449    &  496  & 625    &   779  & 1987  \\
														 & $10^{-6}$ &  474   & 643    &  706  & 819    &   1200 & 2769  \\
														 & $10^{-8}$ &  1017  & 823    &  1674 & 1162   &   2112 & 5000$^*$ \\
\hline
\multirow{3}{*}{$D_k^{XK}$}  & $10^{-4}$ &  112   & 140    &  273  & 240    &   570  & 719   \\
														 & $10^{-6}$ &  160   & 177    &  307  & 312    &   908  & 972   \\
														 & $10^{-8}$ &  317   & 212    &  487  & 332    &   937  & 1172  \\
\hline
\end{tabular}
\end{center}
\end{table}

\begin{figure}[ht]
\begin{center}
\begin{tabular}{cc}
\includegraphics[width=.45\textwidth]{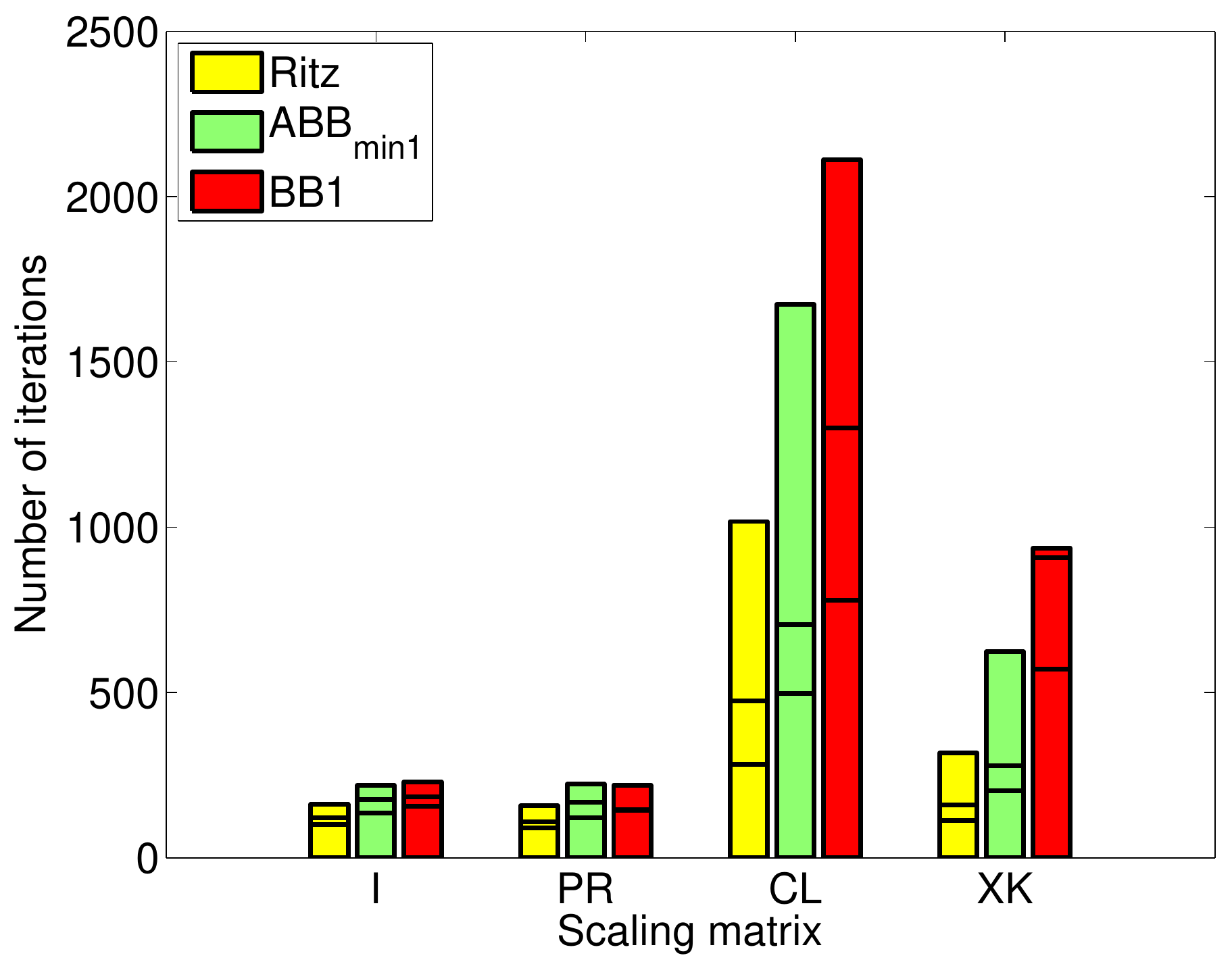}&
\includegraphics[width=.45\textwidth]{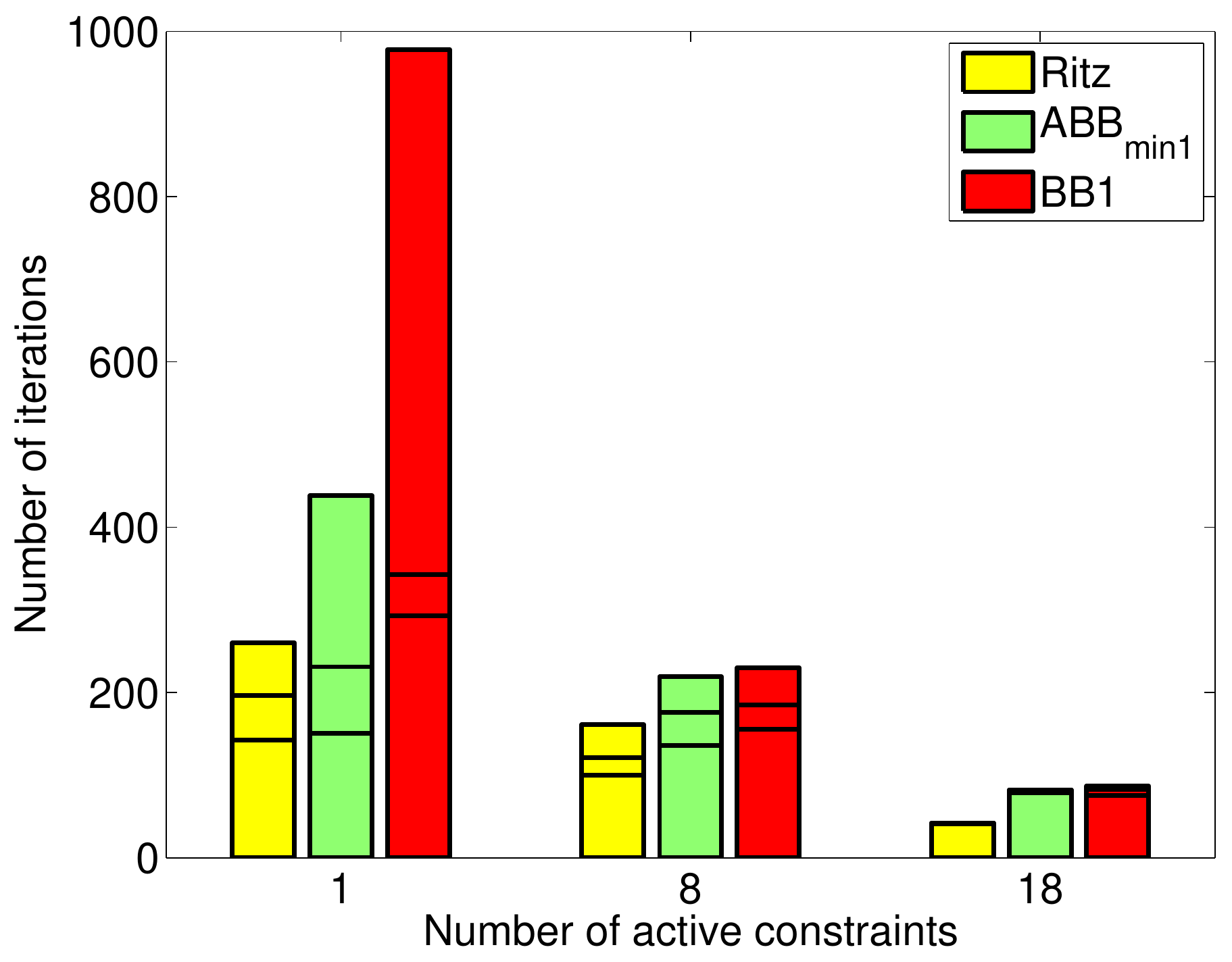}\\
\includegraphics[width=.45\textwidth]{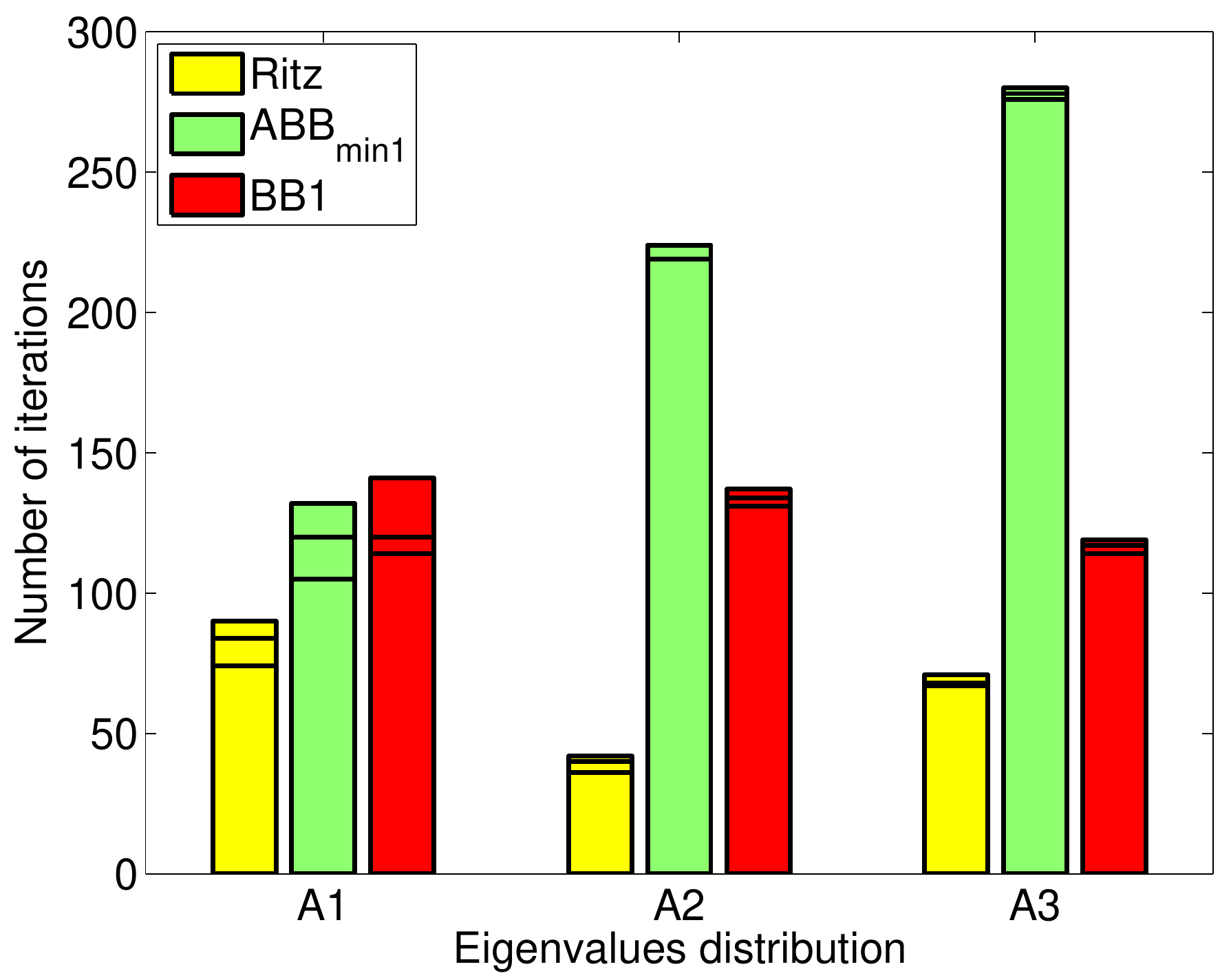}&
\includegraphics[width=.45\textwidth]{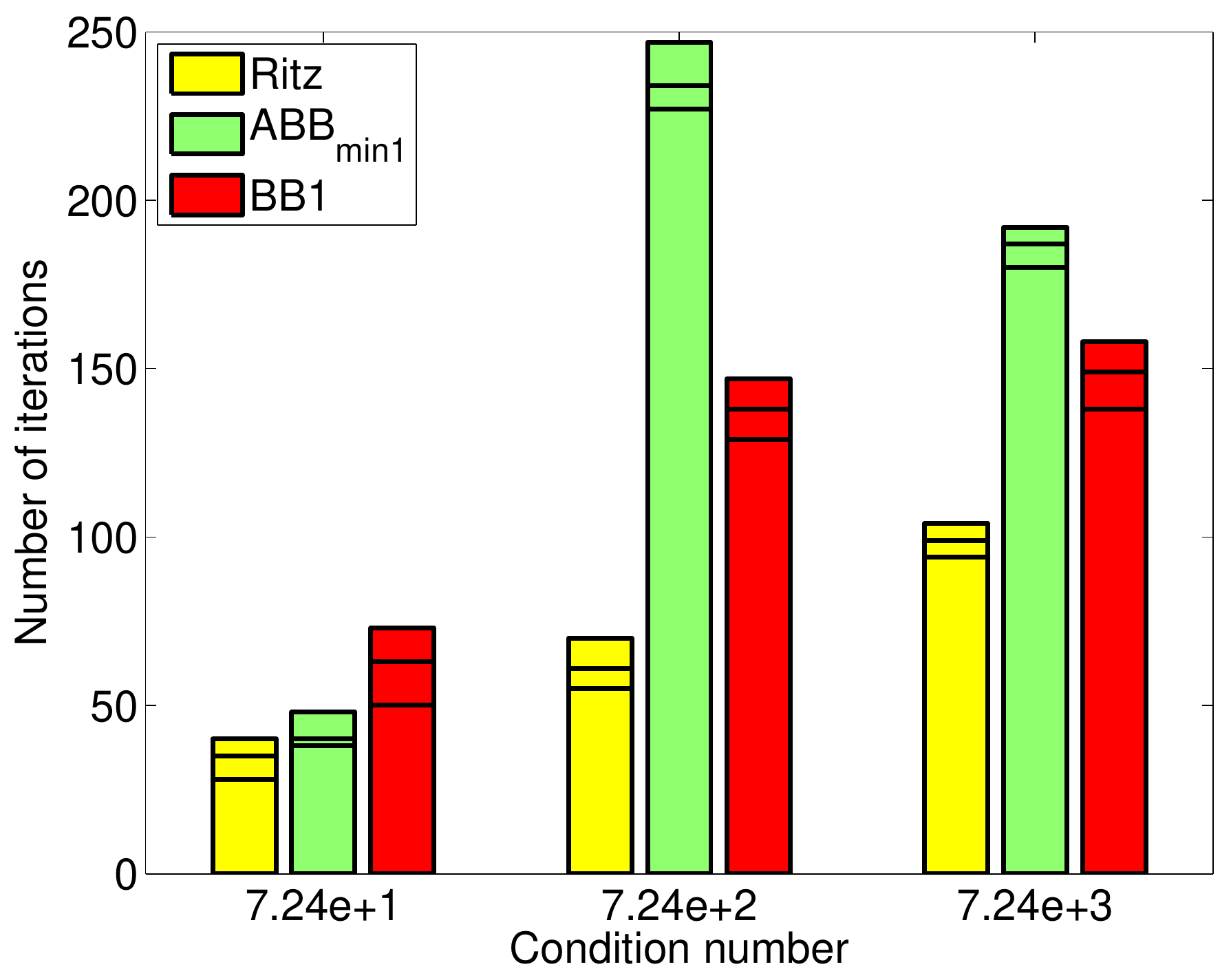}
\end{tabular}
\caption{Graphic representation of the results obtained in the quadratic tests when a monotone linesearch is employed by the algorithms. Top left: variation of the scaling matrix. Top right: variation of the number of active constraints. Bottom left: variation of the eigenvalues distribution. Bottom right: variation of the condition number. The lower (resp. upper) horizontal segment of each bar represents the number of iterations needed to reach a RRE lower than $10^{-4}$ (resp. $10^{-6}$), while the height of each bar corresponds to a RRE $\leq 10^{-8}$.}
\label{figG1}
\end{center}
\end{figure}

\begin{figure}[ht]
\begin{center}
\begin{tabular}{cc}
\includegraphics[width=.45\textwidth]{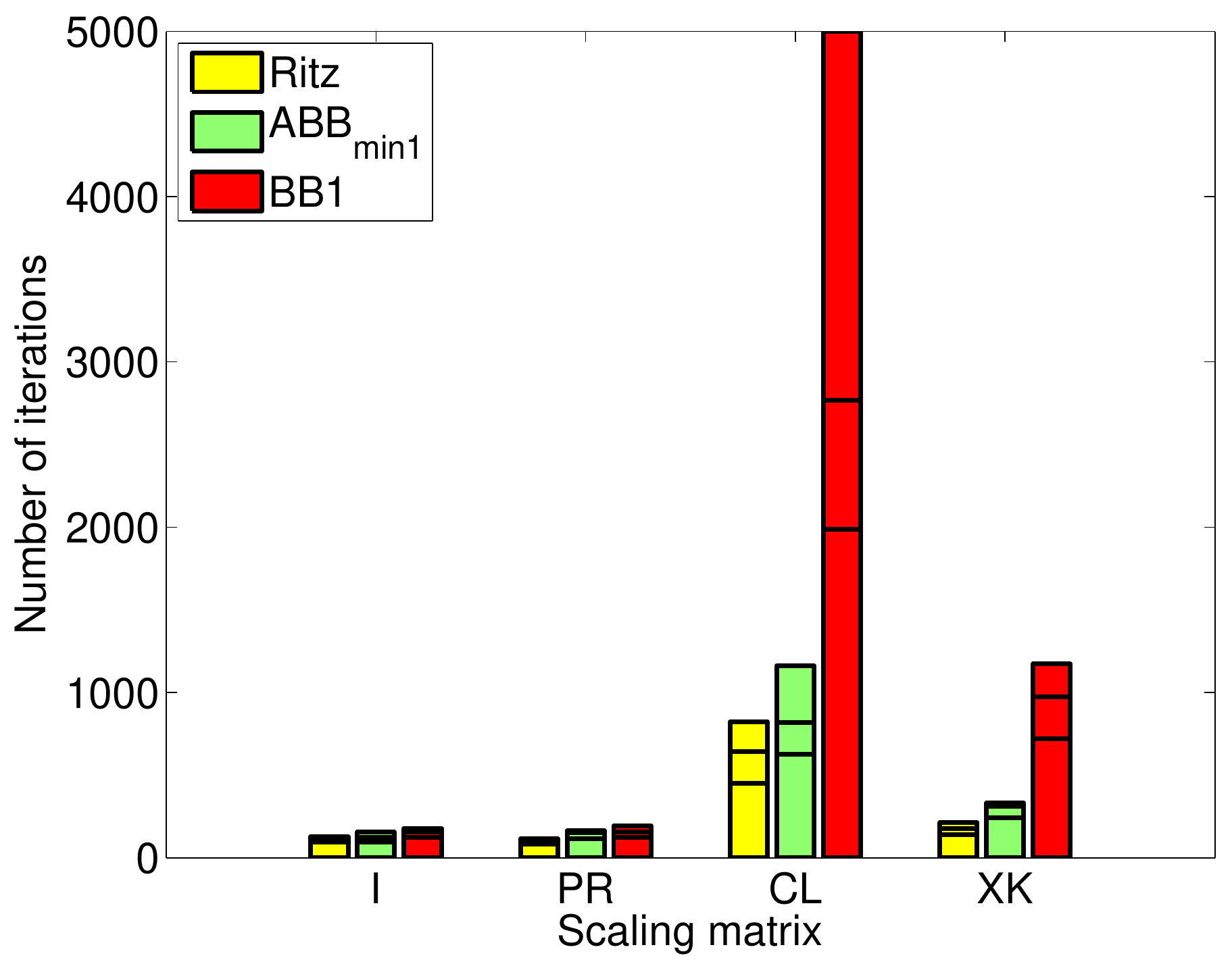}&
\includegraphics[width=.45\textwidth]{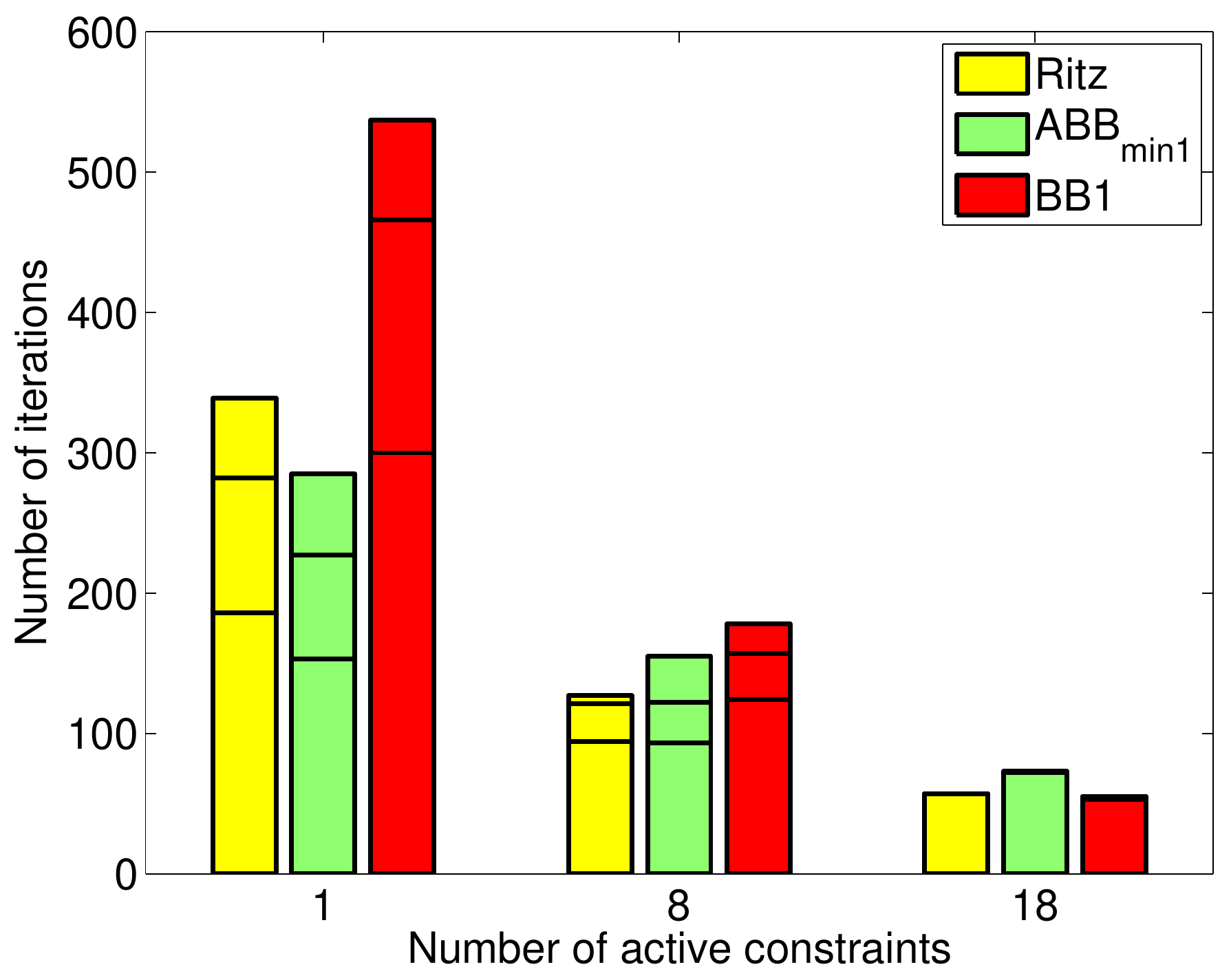}\\
\includegraphics[width=.45\textwidth]{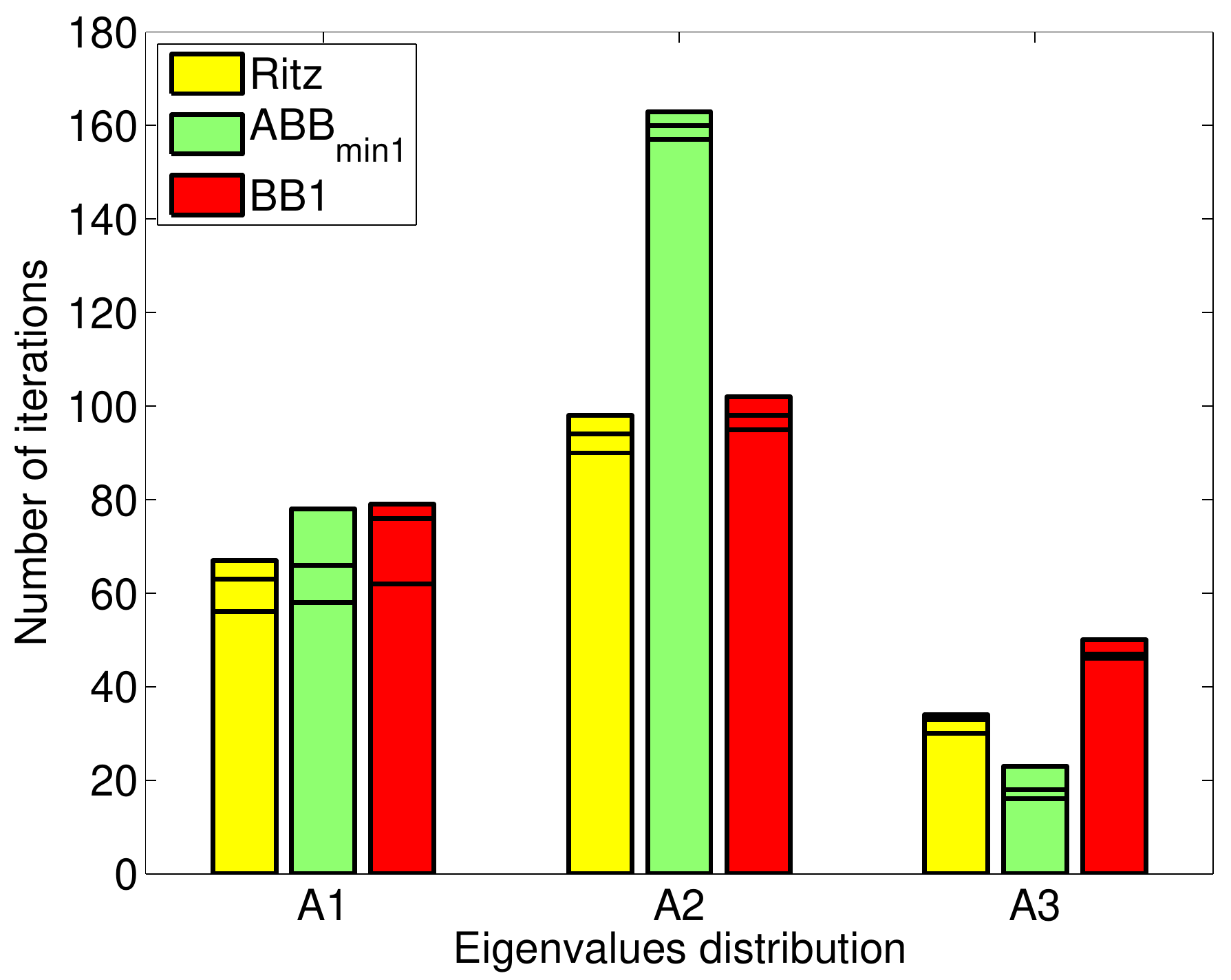}&
\includegraphics[width=.45\textwidth]{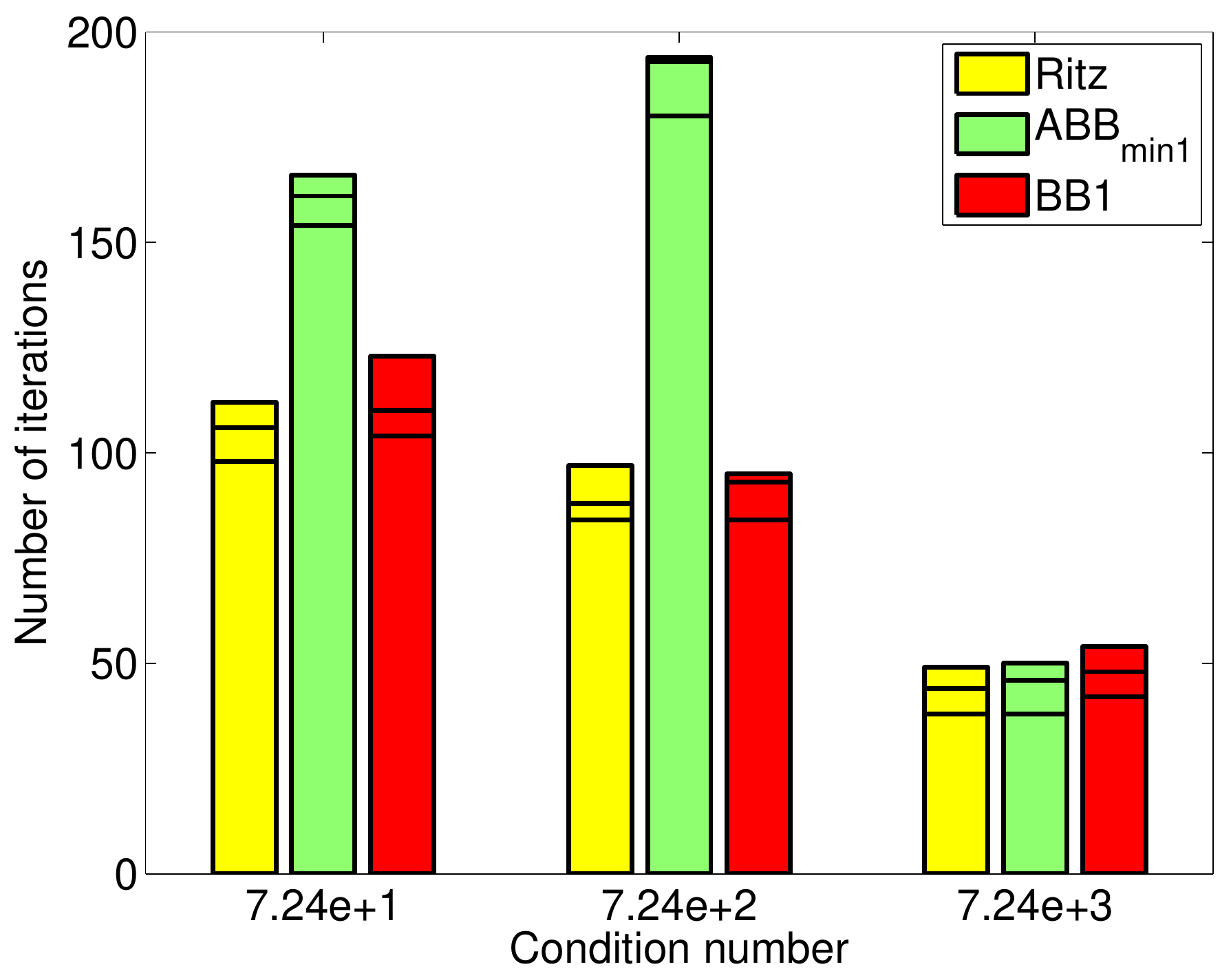}
\end{tabular}
\caption{Graphic representation of the results obtained in the quadratic tests when a nonmonotone linesearch is employed by the algorithms. Top left: variation of the scaling matrix. Top right: variation of the number of active constraints. Bottom left: variation of the eigenvalues distribution. Bottom right: variation of the condition number. The lower (resp. upper) horizontal segment of each bar represents the number of iterations needed to reach a RRE lower than $10^{-4}$ (resp. $10^{-6}$), while the height of each bar corresponds to a RRE $\leq 10^{-8}$.}
\label{figG1b}
\end{center}
\end{figure}

From the information provided in Table \ref{tabG1} and shown graphically in the top left panels of Figures \ref{figG1} and \ref{figG1b} we can see that the choice of the steplength provided by the limited memory rule is able to reduce substantially the number of iterations required to reach a given accuracy, with maximum gains of more than 30\% of iterations with respect to the ABB$_{\rm{min1}}$ strategy. The BB1 steplength seems to be less effective in all cases. As concerns the comparison between the scaling matrices, the best performances are obtained with the two stationary choices (i.e., the identity or the inverse of the diagonal of $A$), while the XK and in particular the CL scaling matrices exhibit a clear slower convergence rate.\\
In the following tests, we used the nonscaled GP algorithm and we analyzed the behaviour of the schemes for:
\begin{itemize}
\item[$\bullet$] different values of the number of active constraints $n_a$: 1, 8, 18. The results are shown in Table \ref{tabG2} and in the top right panels of Figures \ref{figG1} and \ref{figG1b};
\item[$\bullet$] different eigenvalues distributions. We decided to fix the number of active constraints $n_a = 8$ and to keep the condition number unchanged by setting again $\xi_1 = 1$ and $\xi_{20} = (\sqrt{2})^{19}$, and we chose $\xi_2,\ldots,\xi_{19}$ randomly in A1) $(\xi_1,\xi_{20}/3)$, A2) $(\xi_{20}/3,2\xi_{20}/3)$, and A3) $(2\xi_{20}/3,\xi_{20})$, in order to simulate eigenvalues around the minimum value, in the central part of the spectrum and around the maximum value. The results are reported in Table \ref{tabG3} and in the bottom left panels of Figures \ref{figG1} and \ref{figG1b};
\item[$\bullet$] different condition numbers of $A$. To this aim, we fixed again $\xi_{20}=(\sqrt{2})^{19}$ and we changed $\xi_1$ in order to modify the condition number of $A$. Since $(\sqrt{2})^{19} \approx 724$, we tried $\xi_1 = 0.1,1,10$, thus leading to condition numbers of about $7240,724,72.4$, respectively. The eigenvalues $\xi_2,\ldots,\xi_{19}$ have been chosen randomly in $(\xi_1,\xi_{20})$ while we fixed again the number of active constraints of the solution $n_a$ equal to 8. The results are shown in Table \ref{tabG4} and in the bottom right panels of Figures \ref{figG1} and \ref{figG1b}.
\end{itemize}

\begin{table}[ht]
\caption{Numbers of iterations required by GP equipped with the limited memory (Ritz), ABB$_{\rm{min1}}$ and BB1 steplengths to reach RREs lower than $10^{-4}$, $10^{-6}$ and $10^{-8}$ for different numbers of active constraints. The results obtained with a monotone ($M=1$) and nonmonotone ($M=10$) linesearch are reported.}
\label{tabG2}
\begin{center}
\begin{tabular}{c|c|cc|cc|cc|}
\multicolumn{2}{c|}{}           & \multicolumn{2}{c|}{Ritz} & \multicolumn{2}{c|}{ABB$_{\rm{min1}}$} & \multicolumn{2}{c|}{BB1}\\
$n_a$               &    Tol    & M = 1 & M = 10 & M = 1 & M = 10 & M = 1 & M = 10 \\
\hline
\multirow{3}{*}{1}  & $10^{-4}$ &  142  & 186    &  150  & 153    &   293 & 300 \\
		                & $10^{-6}$ &  196  & 282    &  231  & 227    &   343 & 466 \\
		                & $10^{-8}$ &  260  & 339    &  438  & 285    &   978 & 537 \\
\hline
\multirow{3}{*}{8}  & $10^{-4}$ &  100  & 94     &  136  & 93     &   155 & 124 \\
		                & $10^{-6}$ &  121  & 121    &  176  & 122    &   185 & 157 \\
		                & $10^{-8}$ &  161  & 127    &  219  & 155    &   230 & 178 \\
\hline
\multirow{3}{*}{18} & $10^{-4}$ &  41   & 57     &  79   & 72     &   75  & 53  \\
		                & $10^{-6}$ &  42   & 57     &  82   & 72     &   83  & 55  \\
		                & $10^{-8}$ &  42   & 57     &  82   & 73     &   87  & 55  \\
\hline
\end{tabular}
\end{center}
\end{table}

\begin{table}[ht]
\caption{Numbers of iterations required by GP equipped with the limited memory (Ritz), ABB$_{\rm{min1}}$ and BB1 steplengths to reach RREs lower than $10^{-4}$, $10^{-6}$ and $10^{-8}$ for different eigenvalues distributions (see text). The results obtained with a monotone ($M=1$) and nonmonotone ($M=10$) linesearch are reported.}
\label{tabG3}
\begin{center}
\begin{tabular}{c|c|cc|cc|cc|}
\multicolumn{2}{c|}{}            & \multicolumn{2}{c|}{Ritz} & \multicolumn{2}{c|}{ABB$_{\rm{min1}}$} & \multicolumn{2}{c|}{BB1}\\
${\rm{Eig}}$        &    Tol     & M = 1 & M = 10 & M = 1 & M = 10 & M = 1  & M = 10 \\
\hline
\multirow{3}{*}{A1} & $10^{-4}$  &  74   & 56     &  105  & 58     &   114  & 62  \\
		                & $10^{-6}$  &  84   & 63     &  120  & 66     &   120  & 76  \\
	                  & $10^{-8}$  &  90   & 67     &  132  & 78     &   141  & 79  \\
\hline
\multirow{3}{*}{A2} & $10^{-4}$  &  36   & 90     &  219  & 157    &   131  & 95  \\
	                  & $10^{-6}$  &  40   & 94     &  224  & 160    &   134  & 98  \\
	                  & $10^{-8}$  &  42   & 98     &  224  & 163    &   137  & 102 \\
\hline
\multirow{3}{*}{A3} & $10^{-4}$  &  67   & 30     &  276  & 16     &   114  & 46  \\
	                  & $10^{-6}$  &  68   & 33     &  278  & 18     &   117  & 47  \\
	                  & $10^{-8}$  &  71   & 34     &  280  & 23     &   119  & 50  \\
\hline
\end{tabular}
\end{center}
\end{table}

\begin{table}[ht]
\caption{Numbers of iterations required by GP equipped with the limited memory (Ritz), ABB$_{\rm{min1}}$ and BB1 steplengths to reach RREs lower than $10^{-4}$, $10^{-6}$ and $10^{-8}$ for different condition numbers of $A$. The results obtained with a monotone ($M=1$) and nonmonotone ($M=10$) linesearch are reported.}
\label{tabG4}
\begin{center}
\begin{tabular}{c|c|cc|cc|cc|}
\multicolumn{2}{c|}{}              & \multicolumn{2}{c|}{Ritz} & \multicolumn{2}{c|}{ABB$_{\rm{min1}}$} & \multicolumn{2}{c|}{BB1}\\
$K(A)$                & Tol        & M = 1 & M = 10 & M = 1 & M = 10 & M = 1  & M = 10 \\
\hline
\multirow{3}{*}{7240} & $10^{-4}$  &  94   & 98     &  180  & 154    &   138  & 104 \\
                      & $10^{-6}$  &  99   & 106    &  187  & 161    &   149  & 110 \\
                      & $10^{-8}$  &  104  & 112    &  192  & 166    &   158  & 123 \\
\hline
\multirow{3}{*}{724}  & $10^{-4}$  &  55   & 84     &  227  & 180    &   129  &  84 \\
                      & $10^{-6}$  &  61   & 88     &  234  & 193    &   138  &  93 \\
                      & $10^{-8}$  &  70   & 97     &  247  & 194    &   147  &  95 \\
\hline
\multirow{3}{*}{72.4} & $10^{-4}$  &  28   & 38     &  38   & 38     &   50   &  42 \\
                      & $10^{-6}$  &  35   & 44     &  40   & 46     &   63   &  48 \\
                      & $10^{-8}$  &  40   & 49     &  48   & 50     &   73   &  54 \\
\hline
\end{tabular}
\end{center}
\end{table}

The different numerical experiments we carried out lead to similar conclusions. In fact, if the ABB$_{\rm{min1}}$ and BB1 steplengths overtake each other according to the features of the problem, the values provided by the limited memory rule allow a systematic reduction of the iterations required. A further interesting feature that we noticed in all our tests (but we did not reported in the results of the paper for practicality reasons) is that the lower number of iterations required by the proposed rule is always combined with the faster recovery of the active set of the solution.

\subsection{Imaging problems}

\begin{figure}[ht]
\begin{center}
 \begin{tabular}{ccc}
\includegraphics[width=.3\textwidth]{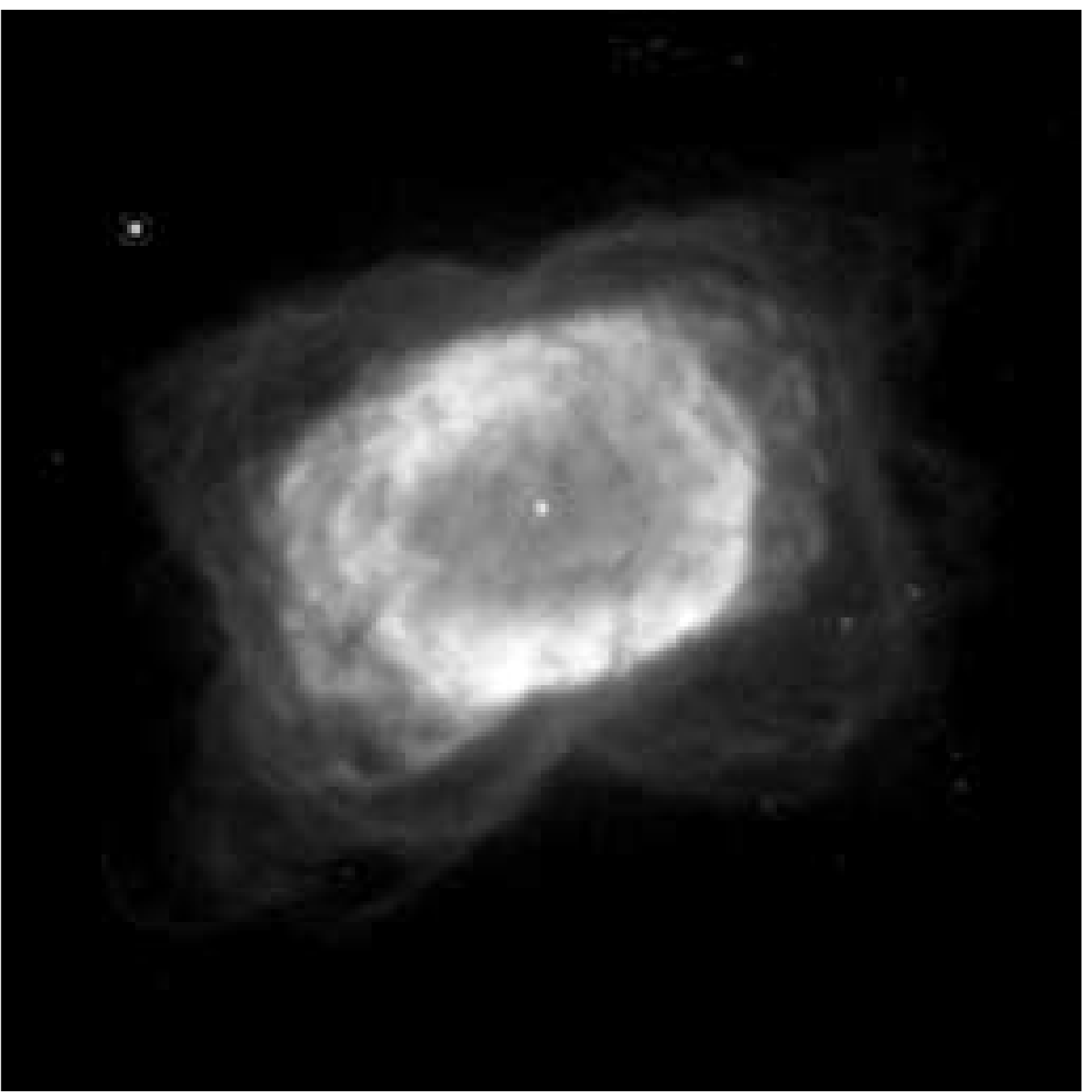}&
\includegraphics[width=.3\textwidth]{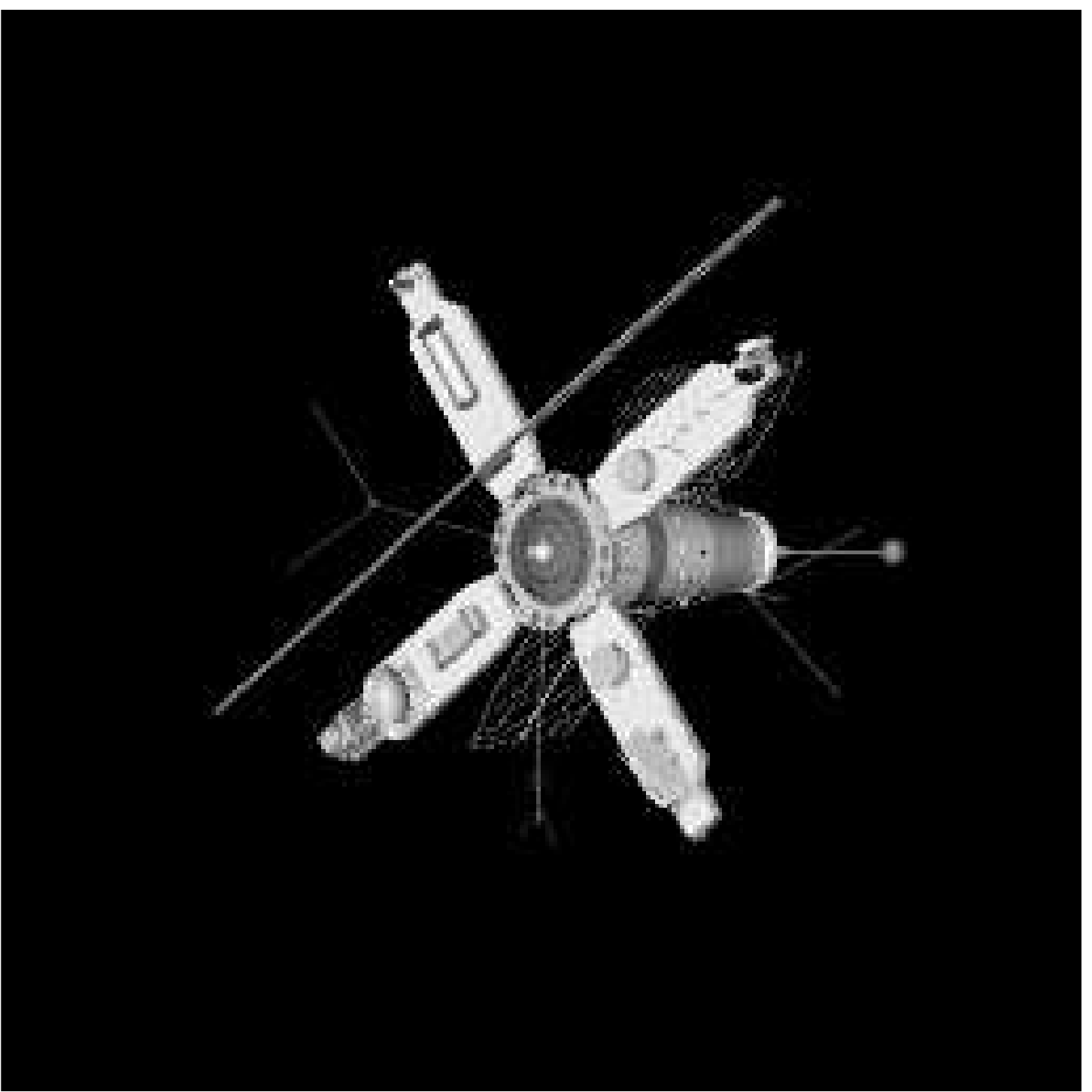}&
\includegraphics[width=.3\textwidth]{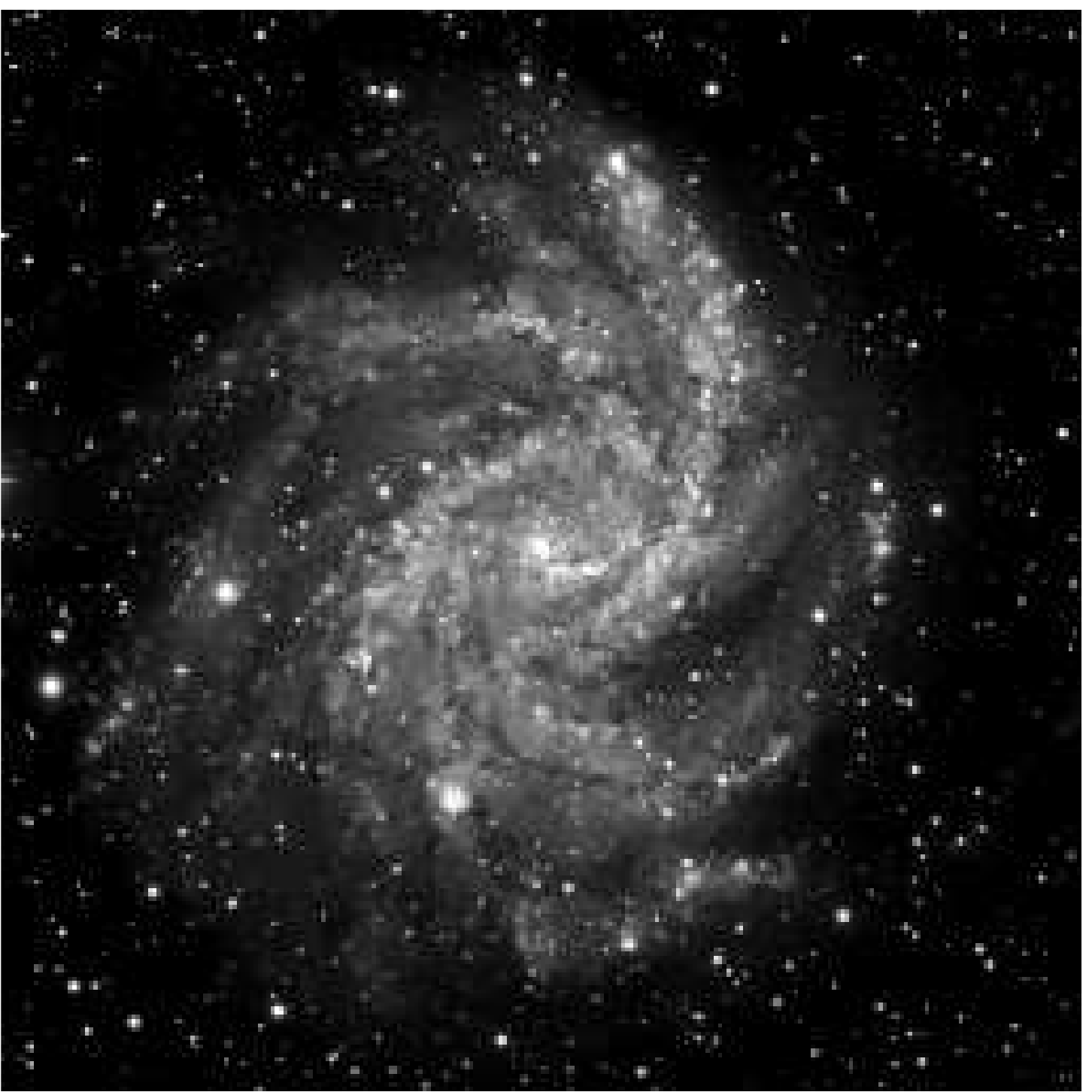}\\
\\
\includegraphics[width=.3\textwidth]{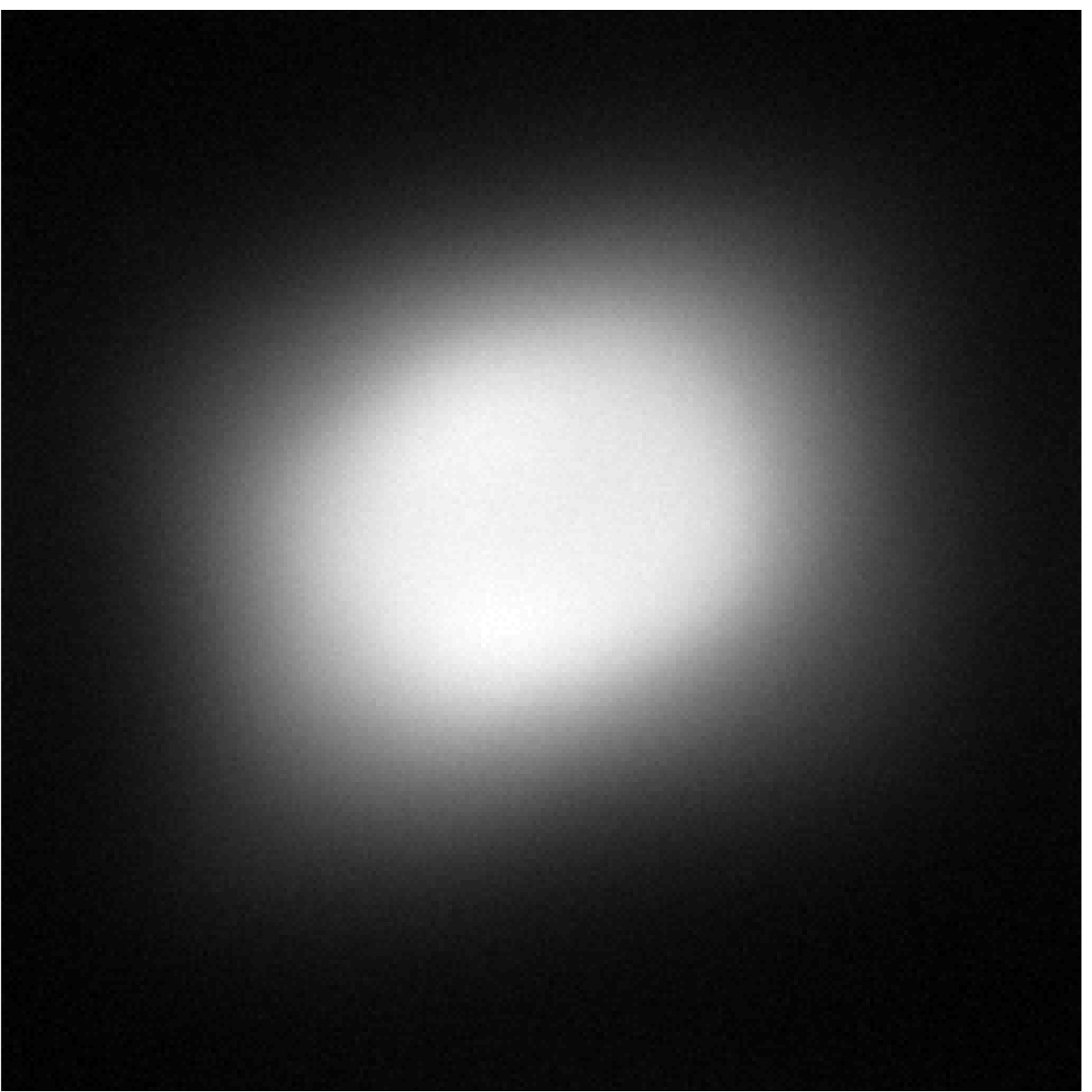}&
\includegraphics[width=.3\textwidth]{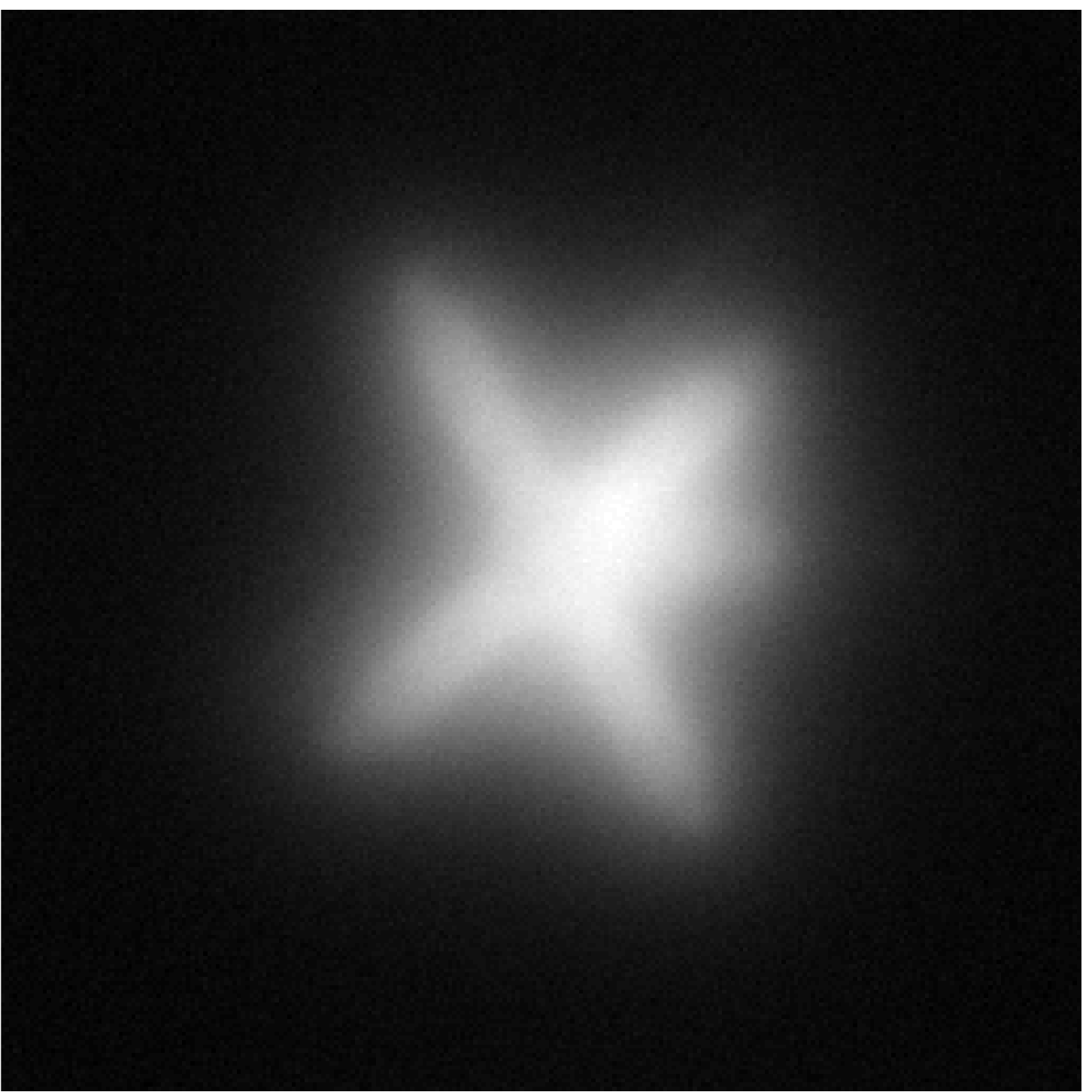}&
\includegraphics[width=.3\textwidth]{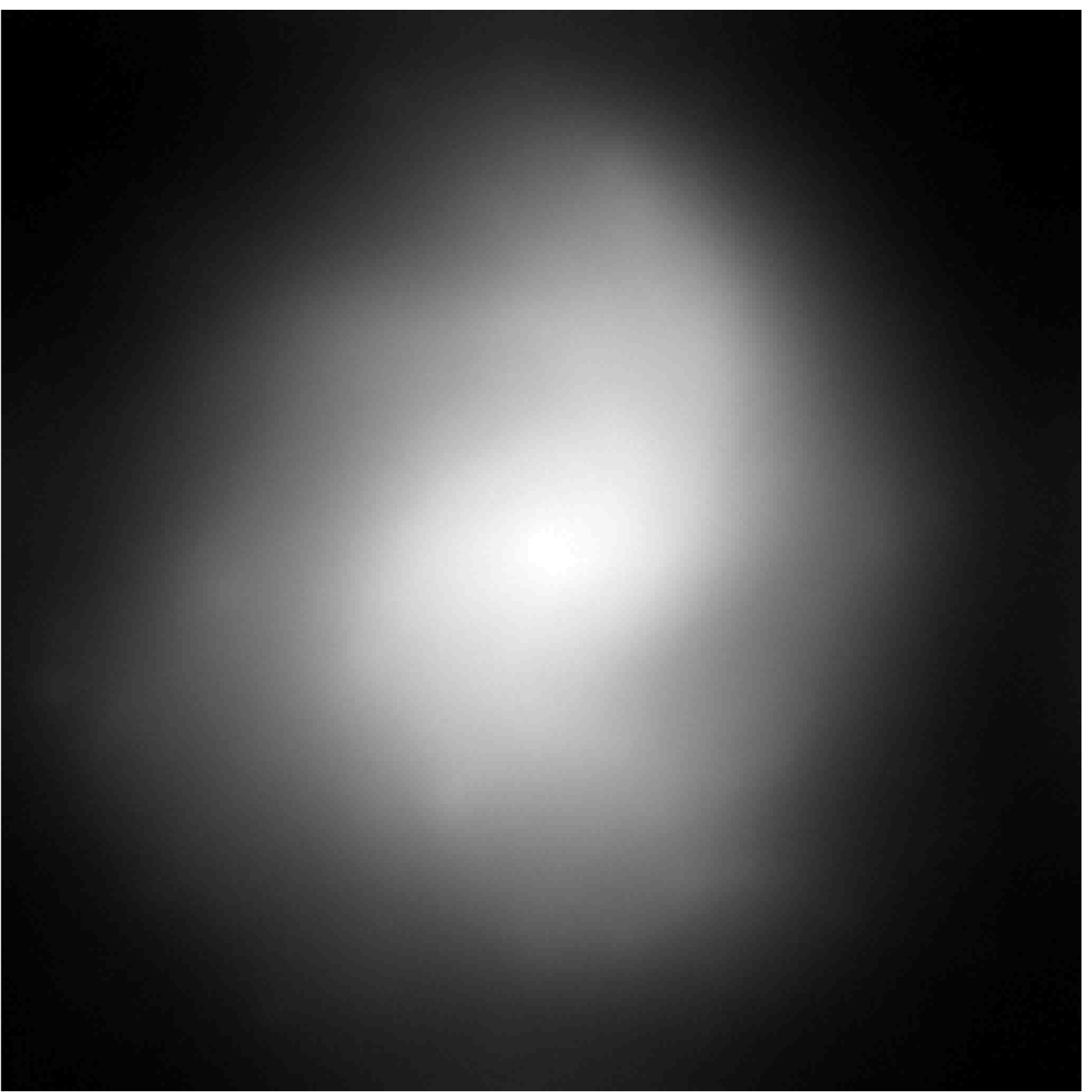}\\
A & B & C\\
 \end{tabular}
\caption{First row: original images for the three test problems. Second row: blurred and noisy images for the three test
problems.}
\label{Immagini}
\end{center}
\end{figure}

In this section we consider a general image reconstruction problem with data perturbed by Gaussian noise and we address the corresponding constrained minimization problem \eqref{min_prob_gen}, with $J \equiv J_0^{LS}$ as defined in \eqref{J_Gaussian}, by means of the following algorithms:
\begin{itemize}
\item[$\bullet$] the nonscaled gradient projection method equipped with either the adaptive BB rule (GP ABB$_{\rm{min1}}$) or the new limited memory steplength selection rule (GP Ritz);
\item[$\bullet$] the scaled gradient projection method equipped with the scaling matrix \eqref{Dk}, with $\beta=0$ and $V_0 \equiv V_0^{LS}$, and either the adaptive BB rule (SGP ABB$_{\rm{min1}}$) or the new limited memory steplength selection rule (SGP Ritz);
\item[$\bullet$] the iterative space reconstruction algorithm (ISRA) \cite{Daube1986}, one of the most exploited method in the literature to deal with the image deblurring problem related to Gaussian noise. ISRA can be seen as a scaled gradient method with constant steplength equal to 1, since its $(k+1)$-th iteration is defined by
\begin{align*}
\ve{x}^{(k+1)} &= {\rm{diag}}\left(\frac{\ve{x}^{(k)}}{A^T(A\ve{x}^{(k)} + \ve{b})}\right)A^T\ve{y} \\
&= \ve{x}^{(k)} - {\rm{diag}}\left(\frac{\ve{x}^{(k)}}{A^T(A\ve{x}^{(k)} + \ve{b})}\right)\nabla J_0^{LS}(\ve{x}^{(k)}).
\end{align*}
\end{itemize}
We point out that, for the (S)GP ABB$_{\rm{min1}}$ approaches, we adopted the modification of the ABB$_{\rm{min1}}$ rule exploited e.g. in \cite{Prato2012,Prato2013}, in which the first 20 steplengths have been chosen equal to the BB2 ones to avoid huge steps at the beginning of the minimization process. Moreover, for all algorithms a monotone linesearch has been adopted to determine the parameter $\lambda_k$. The performances of these methods have been assessed in a comparison with the gradient projection extrapolation (GP Extra) method \cite{Bertsekas2009}, which has the form
\begin{eqnarray}\label{extra_method}
\nonumber \overline{\ve{x}}^{(k)} &=& \ve{x}^{(k)} +\eta_k(\ve{x}^{(k)} - \ve{x}^{(k-1)}),\\
\ve{x}^{(k+1)} &=& \mathbb{P}_+ \left(\overline{\ve{x}}^{(k)} - \alpha \nabla J(\overline{\ve{x}}^{(k)})\right)
\end{eqnarray}
where $\ve{x}^{(-1)} = \ve{x}^{(0)}$ and $\eta_k \in (0,1)$. We will assume that
\begin{equation}\label{eta}
\eta_k = \frac{\theta_k(1-\theta_{k-1})}{\theta_{k-1}}, \qquad k=0,1,...
\end{equation}
where the sequence $\{\theta_k\}$ satisfies $\theta_0 = \theta_1 \in (0,1]$ and
\begin{equation}\label{theta}
\frac{1-\theta_{k+1}}{\theta_{k+1}^2}\leq \frac{1}{\theta_k^2}, \quad \theta_k \leq \frac{2}{k+2}, \quad k=0,1,...
\end{equation}
The following proposition on the iteration complexity of the GP Extra scheme holds true \cite{Bertsekas2009}.
\begin{proposizione}\label{prop}
Let $J:\mathbb{R}^n \longrightarrow \mathbb{R}$ be a convex differentiable function and $d(\ve{x}) = \inf_{\ve{x}^*\in X^*} \| \ve{x} - \ve{x}^* \|, \ve{x} \in \mathbb{R}^n$, where $X^*$ is the set of minimizers of $J$ over the feasible set. Assume that $\nabla J$ is Lipschitz continuous with Lipschitz constant $L$, $X^*$ is nonempty and $J^*$ is the minimum of $J$. Let $\{\ve{x}^{(k)}\}$ be a sequence generated by the algorithm \eqref{extra_method}, where $\alpha = \frac{1}{L}$ and $\eta_k$ satisfies Eqs. \eqref{eta}-\eqref{theta}. Then $\lim_{k\rightarrow \infty} d(\ve{x}^{(k)}) = 0$ and
\begin{equation*}
J(\ve{x}^{(k)}) - J^* \leq \frac{2L}{(k+1)^2}d(\ve{x}^{(0)})^2, \qquad k=1,2,...
\end{equation*}
\end{proposizione}
We remark that the function $J_0^{LS}$ defined in \eqref{J_Gaussian} satisfies the condition required by the GP Extra algorithm for the convergence. In particular, the (smallest) Lipschitz constant of the gradient $\nabla J_0^{LS}$ is
\begin{equation}\label{Lip_Gauss}
L(J_0^{LS}) = \xi_{max}(A^TA),
\end{equation}
where $\xi_{max}(X)$ indicates the maximum eigenvalue of $X$.\\
The test problems here considered are generated by convolving the original $256\times 256$ images, shown in the first row of Figure \ref{Immagini} and denoted by A, B, C, with a point spread function (PSF) and perturbing the results with additive white Gaussian noise with variance 1 (we assume that no background radiation is present). The PSF we adopted is a simulation of a ground-based telescope and can be downloaded from the website http://www.mathcs.emory.edu/$\sim$nagy/ RestoreTools/index.html. For each of the considered images, we show the blurred and noisy data used in the experiments in the second row of Figure \ref{Immagini}.\\

\begin{figure}[ht]
\begin{center}
\begin{tabular}{cc}
\includegraphics[width=.45\textwidth]{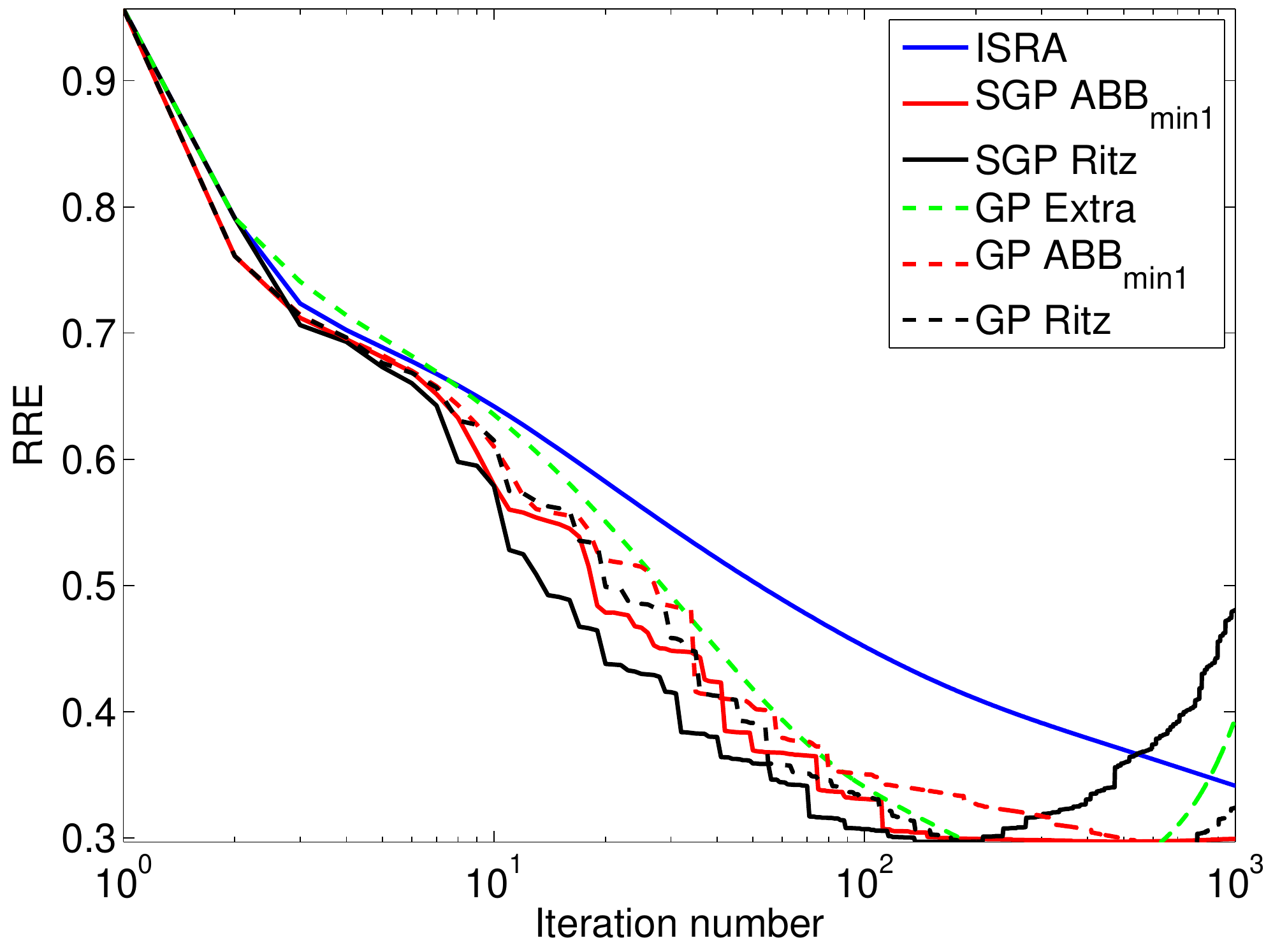}&
\includegraphics[width=.45\textwidth]{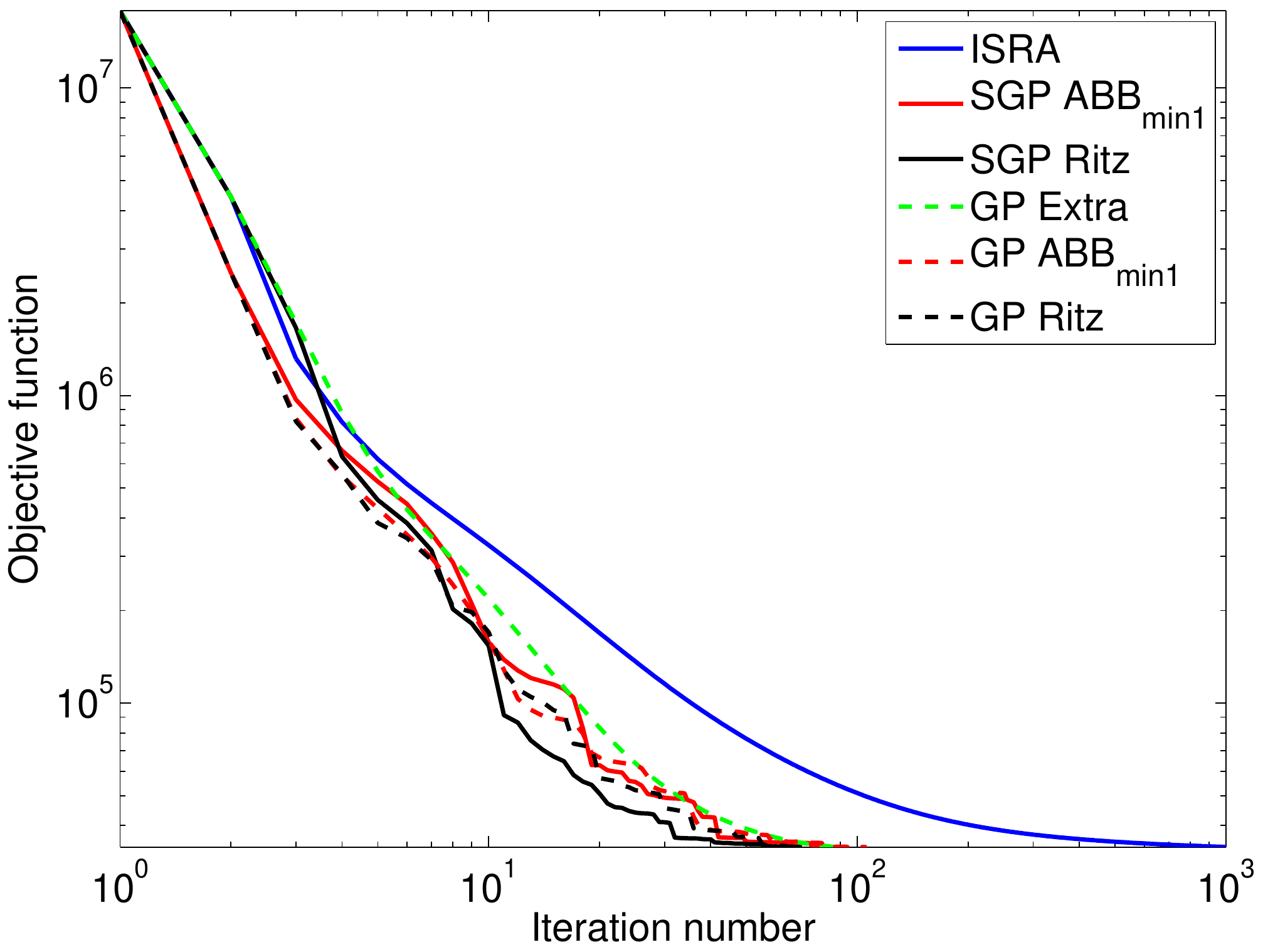}\\
\end{tabular}
\caption{Relative reconstruction error (left panel) and objective function (right panel) provided by the different methods for Image B test problem.}
\label{Err_Funct}
\end{center}
\end{figure}

\begin{figure}[ht]
\begin{center}
\begin{tabular}{ccc}
\includegraphics[width=.3\textwidth]{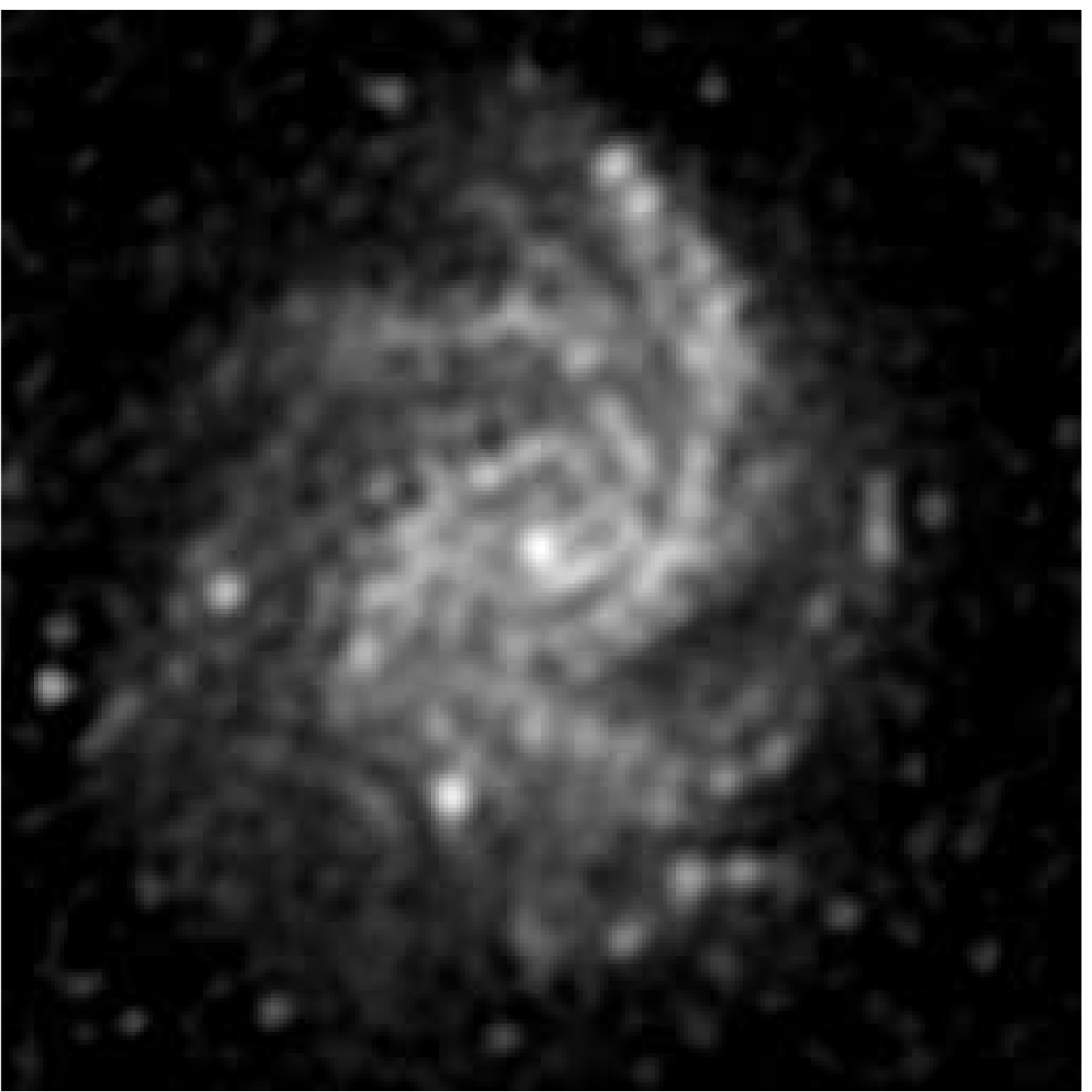}&
\includegraphics[width=.3\textwidth]{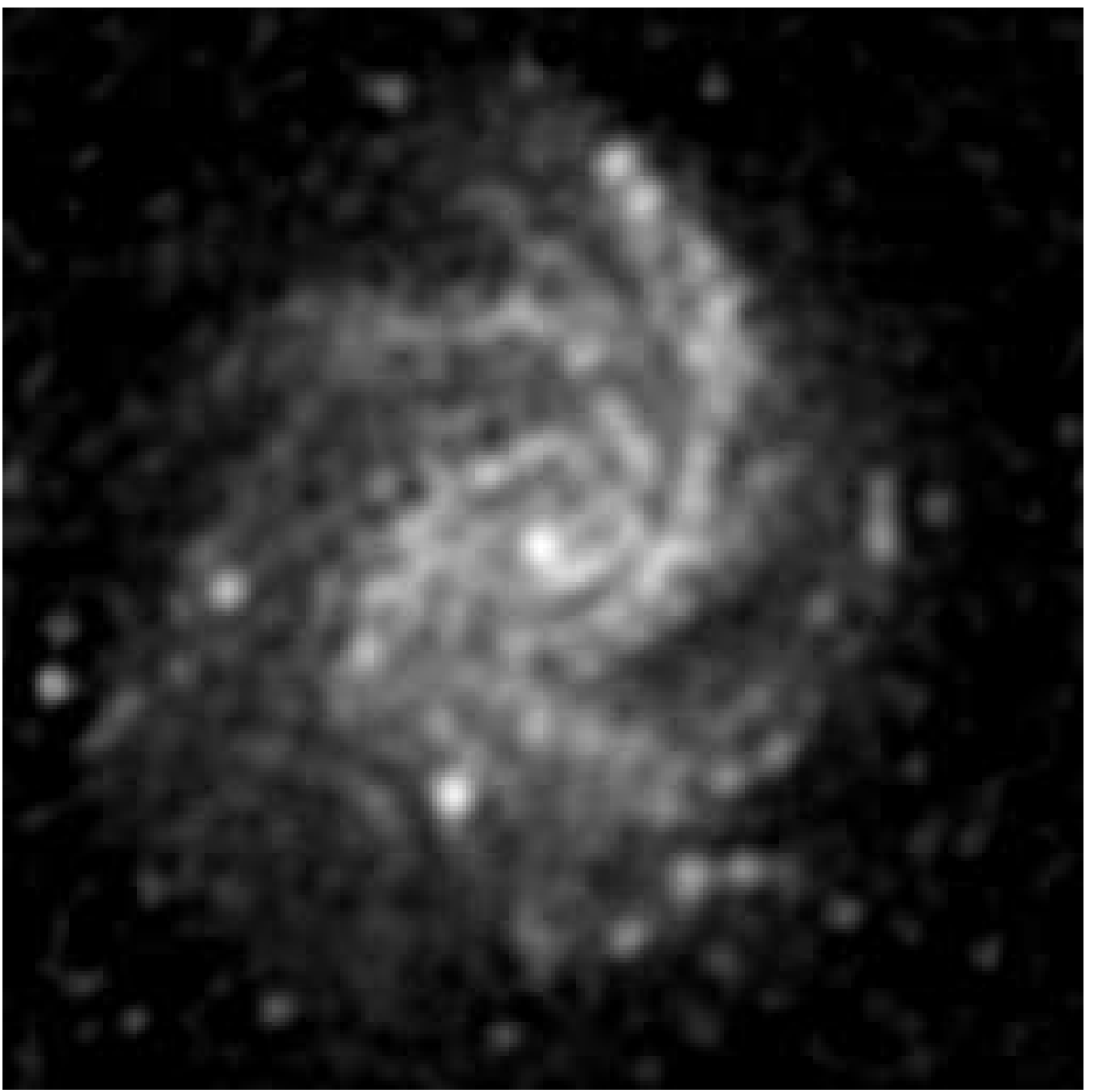}&
\includegraphics[width=.3\textwidth]{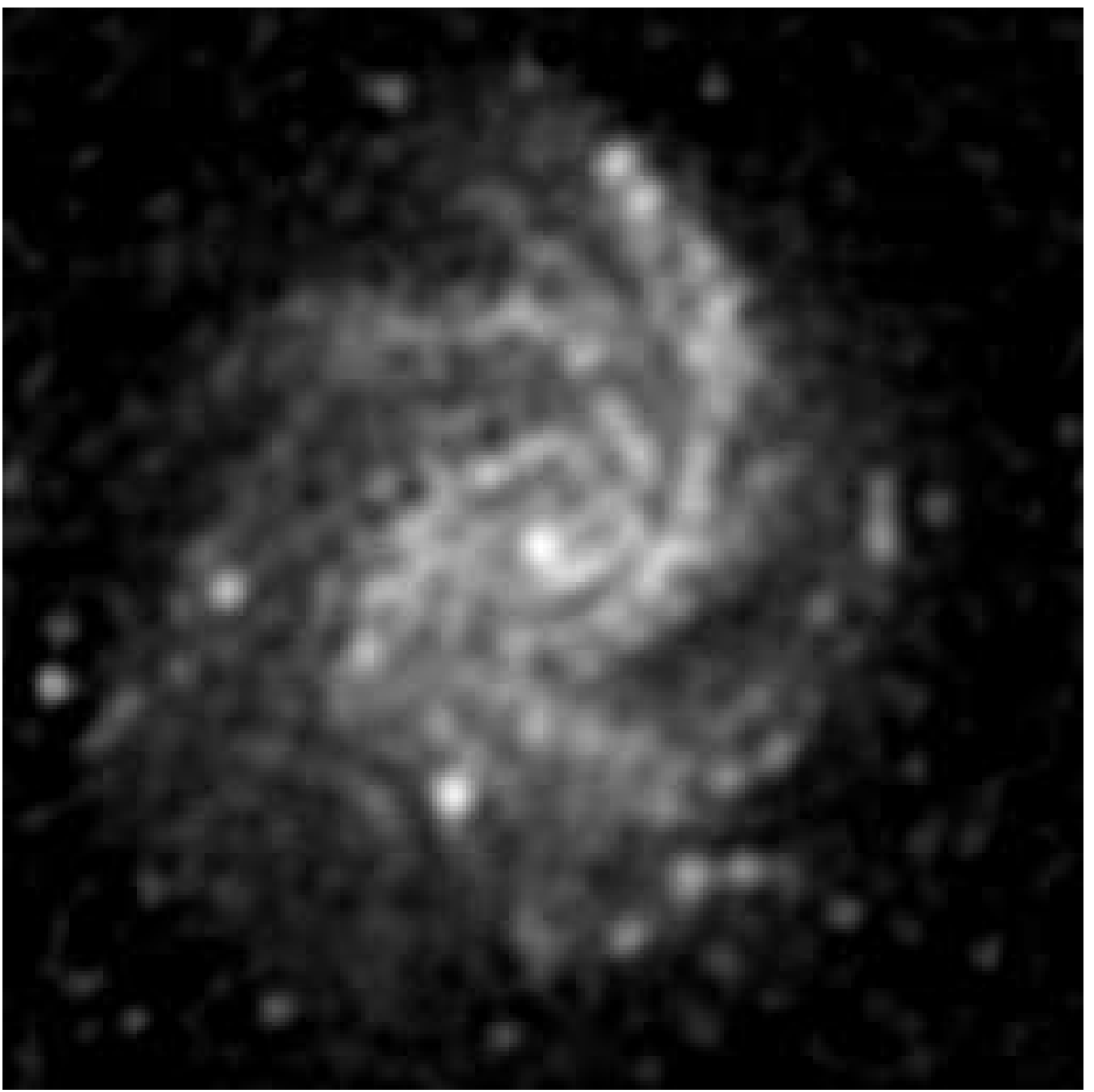}\\
GP Extra & GP ABB$_{\rm{min1}}$ & GP Ritz \\
\includegraphics[width=.3\textwidth]{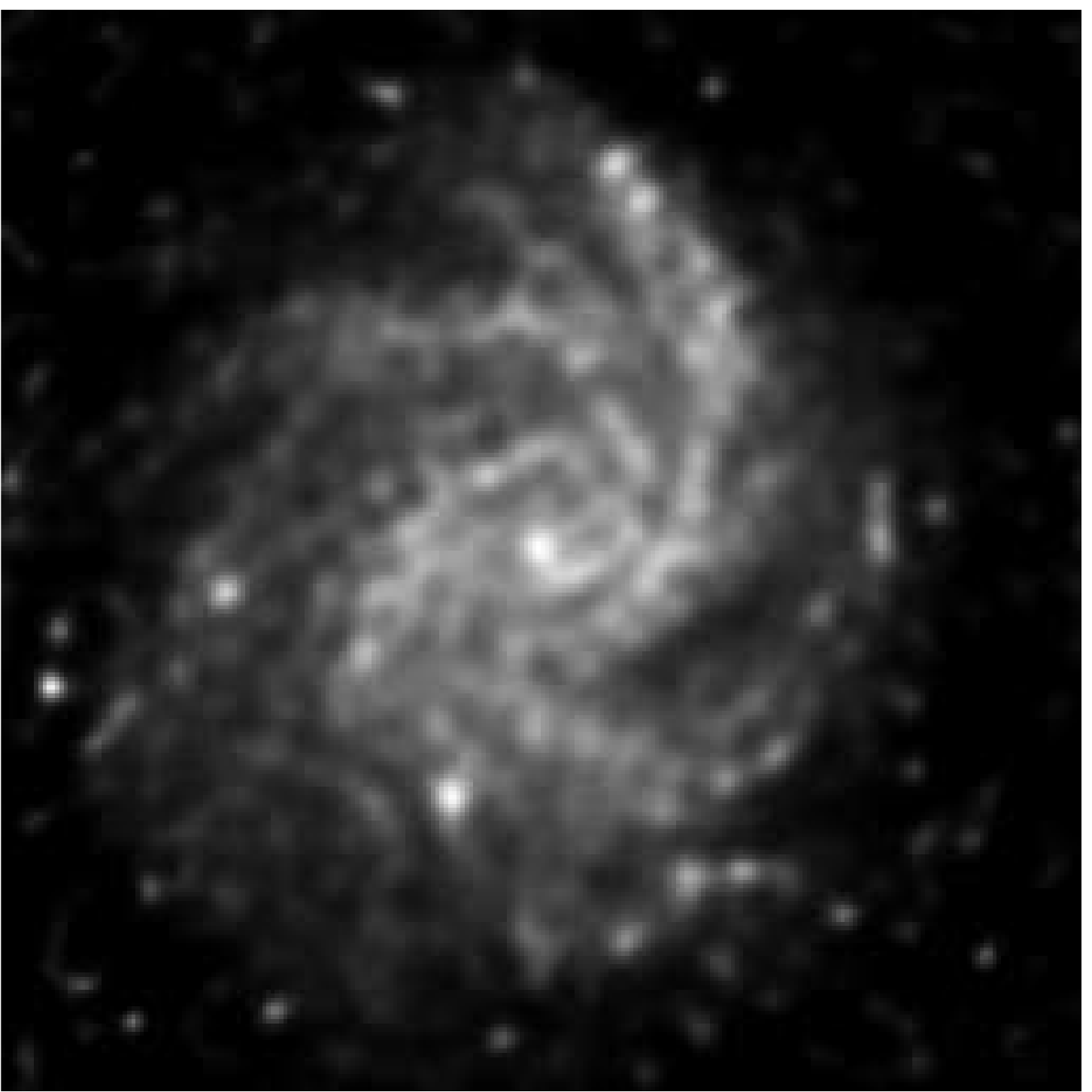}&
\includegraphics[width=.3\textwidth]{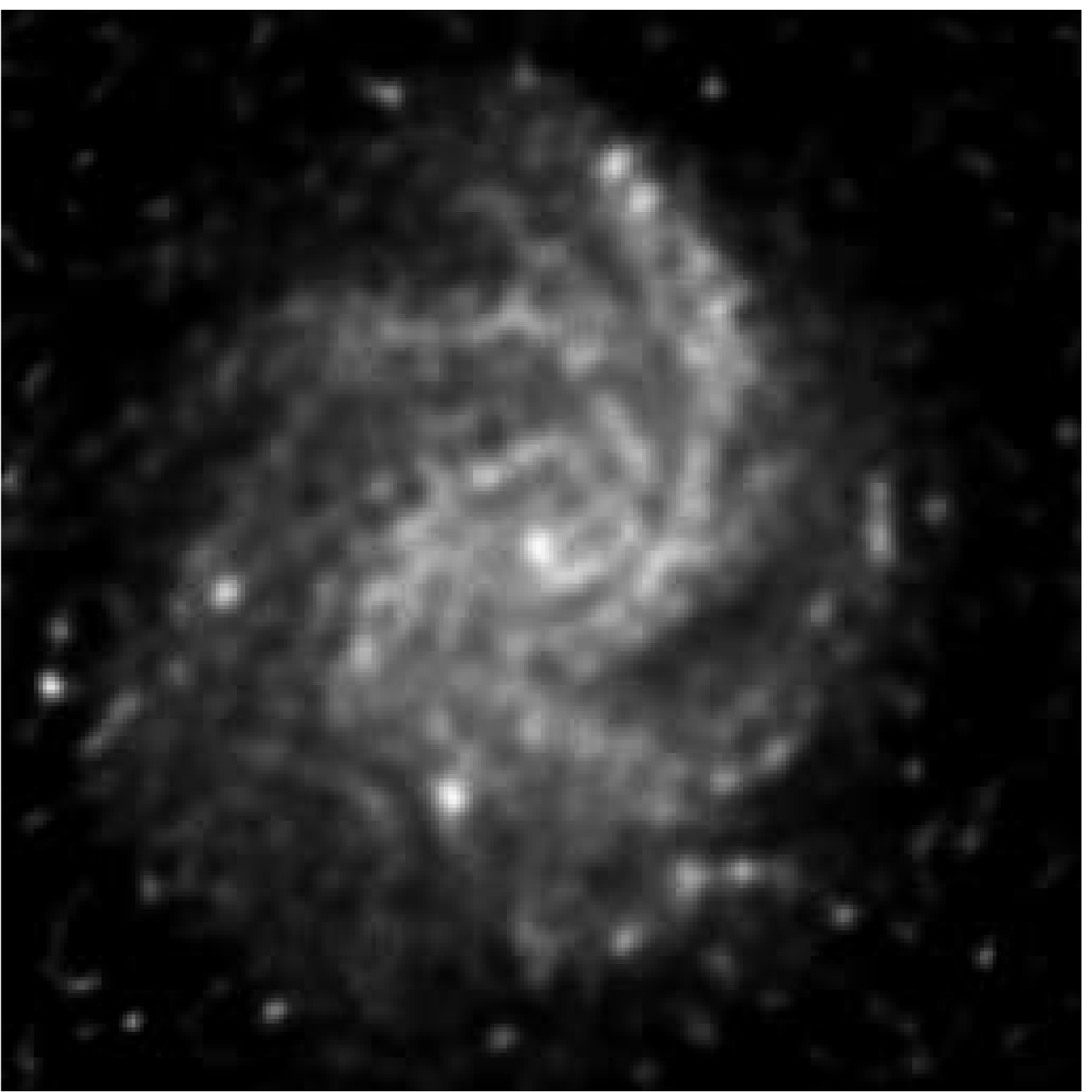}&
\includegraphics[width=.3\textwidth]{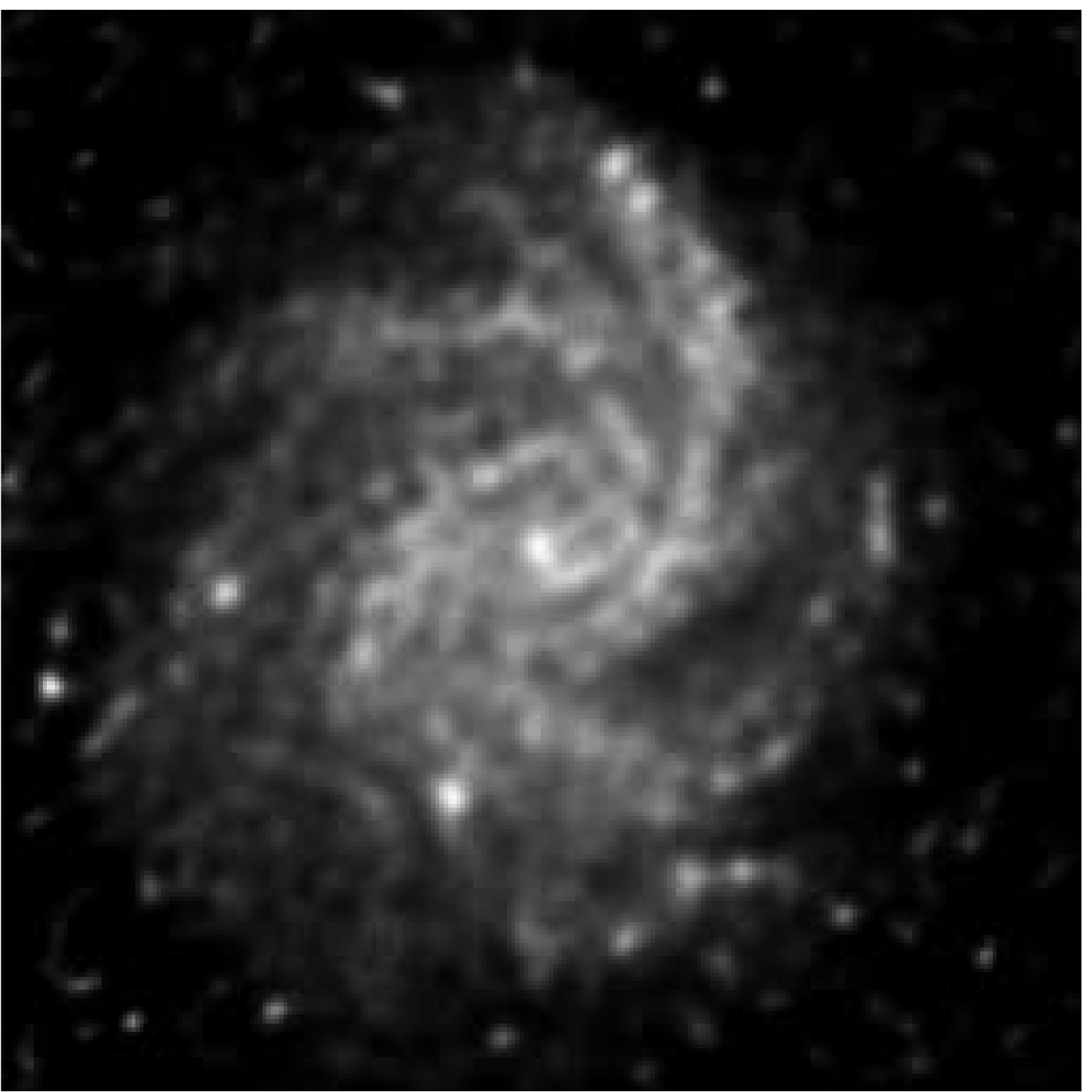}\\
ISRA & SGP ABB$_{\rm{min1}}$ & SGP Ritz \\
\end{tabular}
\caption{Reconstruction of the object C obtained by different GP and SGP methods.}
\label{Rec_Gauss}
\end{center}
\end{figure}

Table \ref{Tab1} reports the minimum RREs and the numbers of iterations required to provide the minimum error, together with the execution times. We remark that, for the GP Extra, the steplength $\alpha$ has been chosen as the reciprocal of the value suggested in \eqref{Lip_Gauss} at each iteration. We show in Figure \ref{Err_Funct} the RRE (as a function of the number of iterations) and the decrease of the objective function provided the different methods for the test problem B, but we appreciated an analogous behavior also from the analysis of problems A and C. To illustrate an example of reconstruction quality provided by a gradient projection method we show the recovered images for test problem C in Figure \ref{Rec_Gauss}.\\
The experiments carried out in section \ref{secG1} are intrinsically different from the tests on imaging problems, since here we do not ask a method to approximate as fast as possible the solution of the minimum problem, but we look for methods which, starting from a given $\ve{x}^{(0)}$, are able to generate a route toward a minimizer of the objective function which passes as close as possible to the original image. The main difference between the two sets of results we obtained is that, in image reconstruction problems, the presence of the scaling matrix has a positive effect (see also \cite{Bonettini2013,Cornelio2013}), as attested by both the lower RREs and the reduced number of iterations needed by the SGP ABB$_{\rm{min1}}$ and SGP Ritz methods with respect to their nonscaled versions GP ABB$_{\rm{min1}}$ and GP Ritz. The constant behaviour noticed on the several tests is that the iterations required by the steplength defined by the limited memory approach are again fewer than those of the alternated scheme, and comparable with a state-of-the-art method as the GP Extra algorithm. It is worth noting that the decrease of the objective function exhibited by the GP and SGP approaches is very similar to that of the GP Extra method, whose iteration complexity has been proved to be $O(1/\sqrt{\varepsilon})$ (see Proposition \ref{prop}). Nevertheless, besides the product $A^TA\ve{x}^{(k)}$ which has to be computed at each iteration by all the algorithms, we have to remark that the GP Extra does not require any additional vector-vector product, with a result of a faster execution time even in cases in which a higher number of iterations are required to provide the best reconstruction (see Table \ref{Tab1}, problems A and C).

\begin{table}[ht]
\caption{Minimum RRE achieved by each algorithm in the Gaussian deblurring problems, with the corresponding number of iterations required and execution time. The asterisk denotes the maximum number of iterations allowed.}
\label{Tab1}
\begin{center}
\setlength{\tabcolsep}{3.5pt}
\begin{tabular}{l|ccc|ccc|ccc|}
                       & \multicolumn{3}{c|}{{\bf A}} & \multicolumn{3}{c|}{{\bf B}} & \multicolumn{3}{c|}{{\bf C}}\\
											 & It. & RRE & Time(s) & It. & RRE & Time(s) & It. & RRE & Time(s)\\
\hline
GP Extra               & 120  & 0.080 & 0.541 &  384 & 0.274 & 1.477 &   420   & 0.296 & 1.665 \\
GP ABB$_{\rm{min1}}$   & 120  & 0.080 & 0.755 & 1684 & 0.276 & 8.671 &  1684   & 0.296 & 9.371 \\
GP Ritz                & 124  & 0.080 & 0.837 &  427 & 0.277 & 2.538 &   839   & 0.296 & 4.896 \\
ISRA                   & 1904 & 0.074 & 5.696 & 5954 & 0.299 & 17.63 & 10000$^*$ & 0.291 & 29.01 \\
SGP ABB$_{\rm{min1}}$  & 90   & 0.074 & 0.753 &  461 & 0.297 & 3.416 &  1003   & 0.288 & 6.726 \\
SGP Ritz               & 91   & 0.074 & 0.702 &  164 & 0.297 & 1.249 &   380   & 0.288 & 2.562 \\
\hline
\end{tabular}
\end{center}
\end{table}

\section{Numerical experiments - Poisson noise}\label{sec5}

For the case of image reconstruction problems with data affected by Poisson noise, we evaluated the utility of the limited memory steplength selection rule in both the presence and the absence of an explicit regularization term in the objective function.

\subsection{Approach without regularization terms}\label{not_reg_sec}

\begin{figure}[htb!]
\begin{center}
\begin{tabular}{cc}
\includegraphics[width=.3\textwidth]{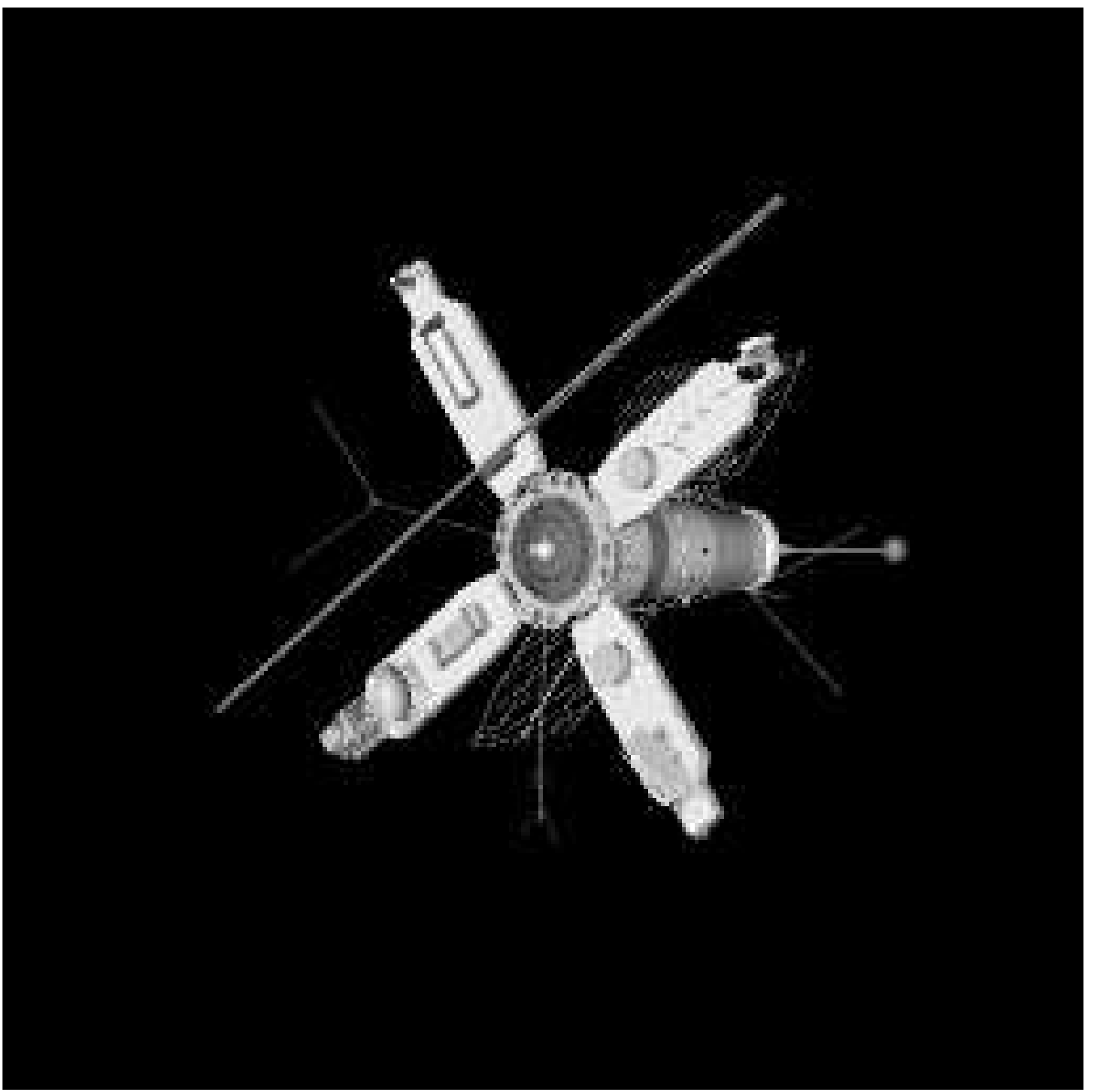}&
\includegraphics[width=.3\textwidth]{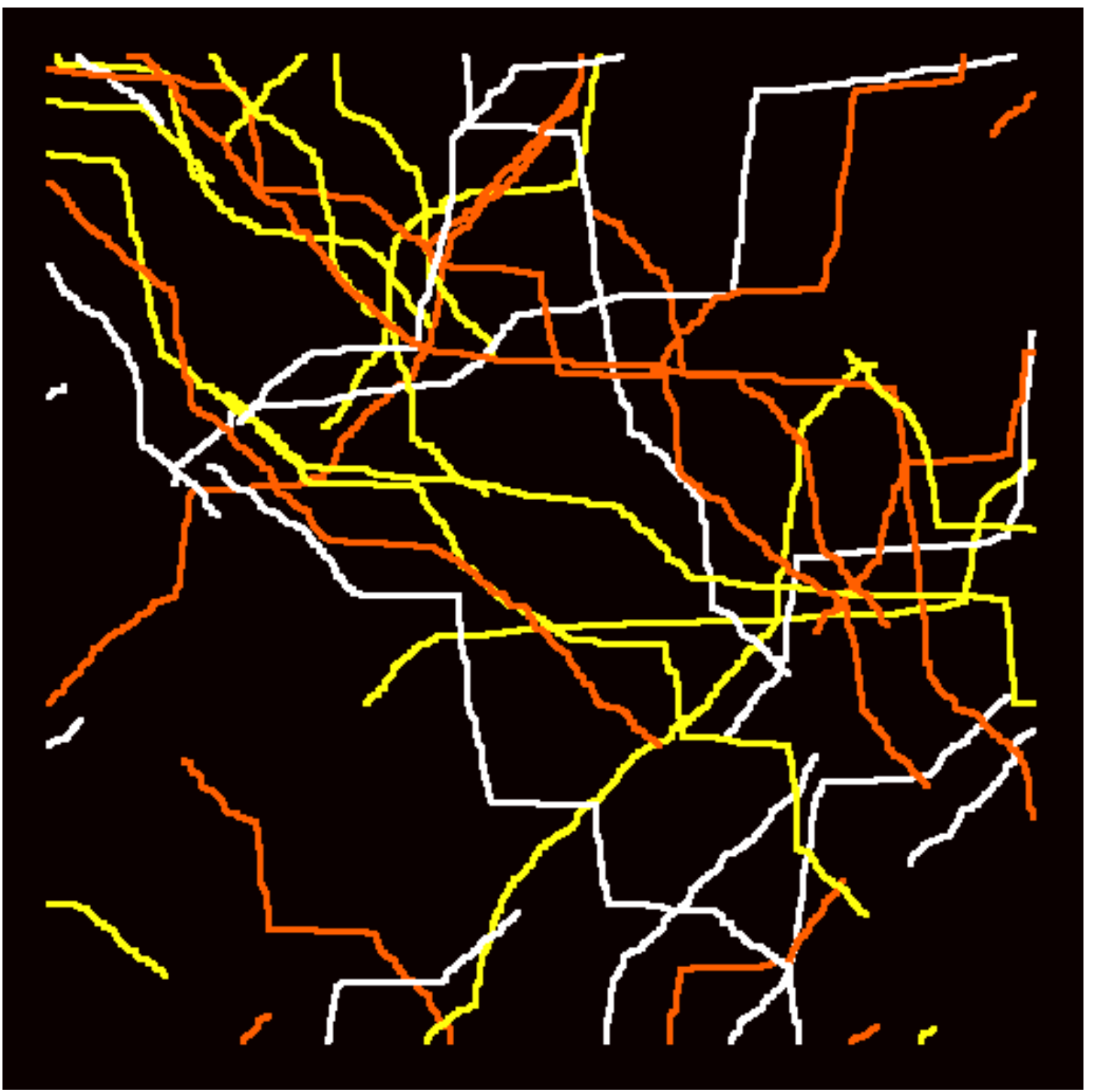}\\
\includegraphics[width=.3\textwidth]{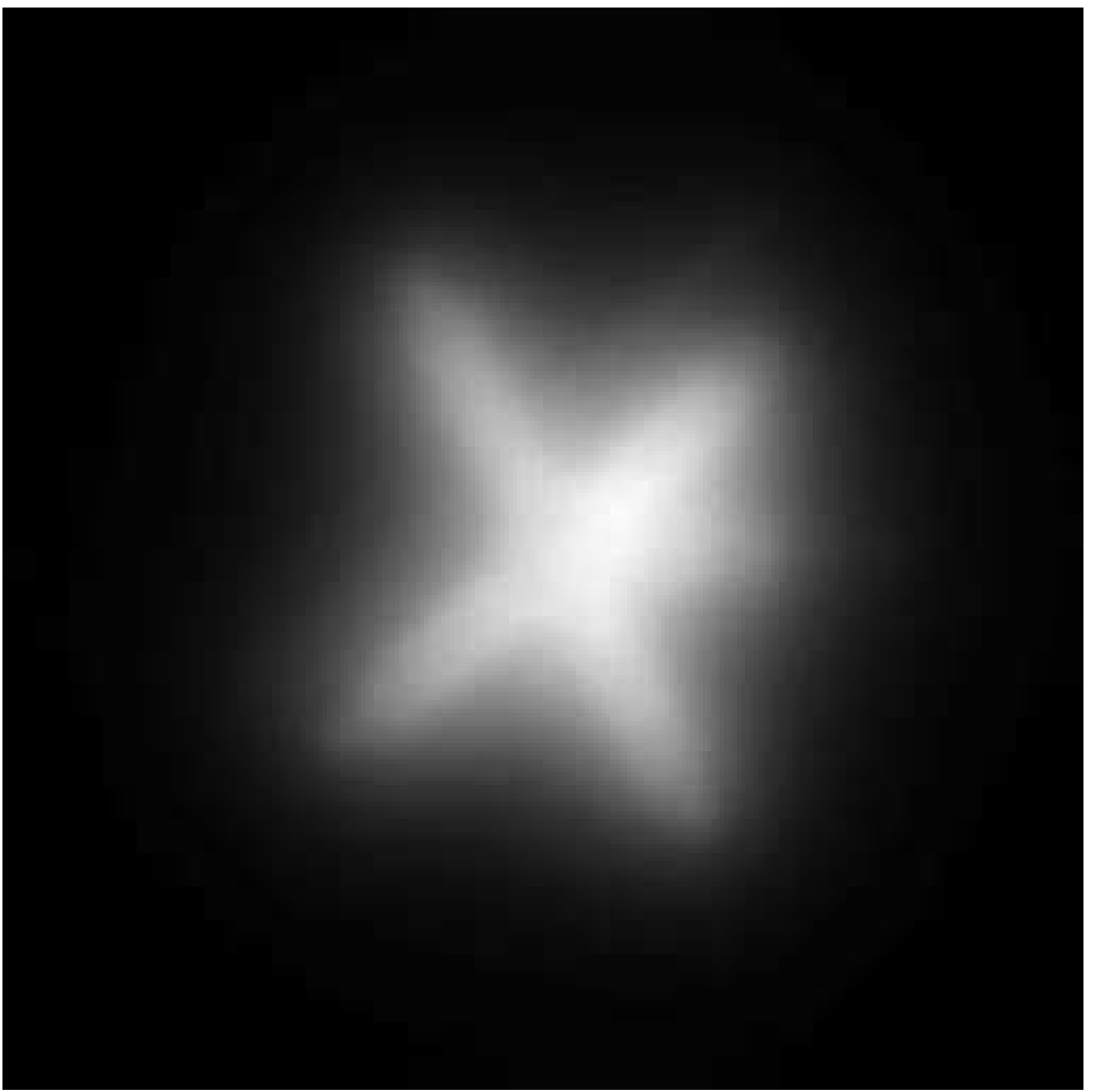}&
\includegraphics[width=.3\textwidth]{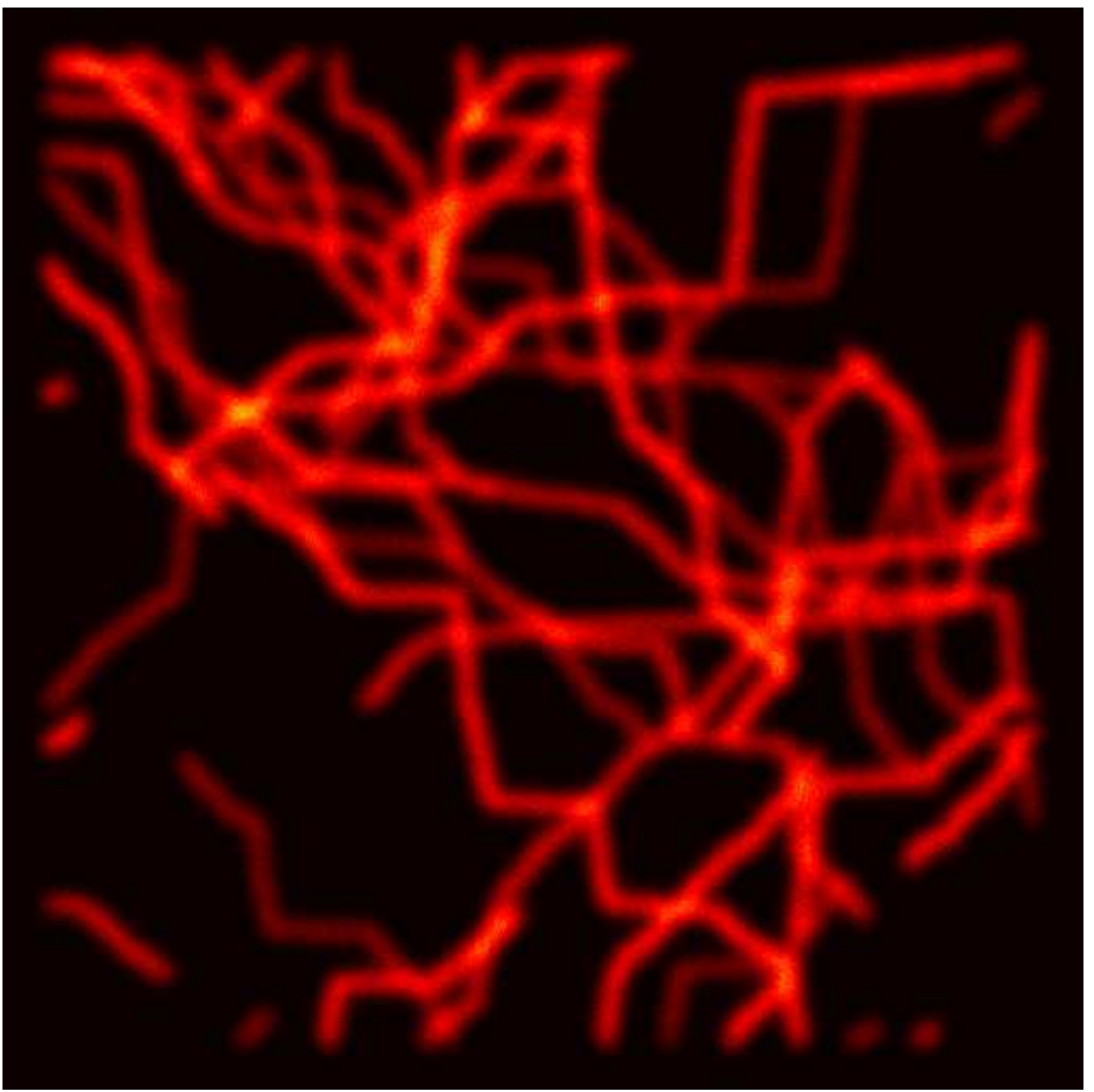}\\
D & E\\
\end{tabular}
\caption{First row: original images for the two test problems. Second row: blurred and noisy images for the two test problems.}
\label{Immagini2}
\end{center}
\end{figure}

\begin{figure}[ht]
\begin{center}
\begin{tabular}{cc}
\includegraphics[width=.45\textwidth]{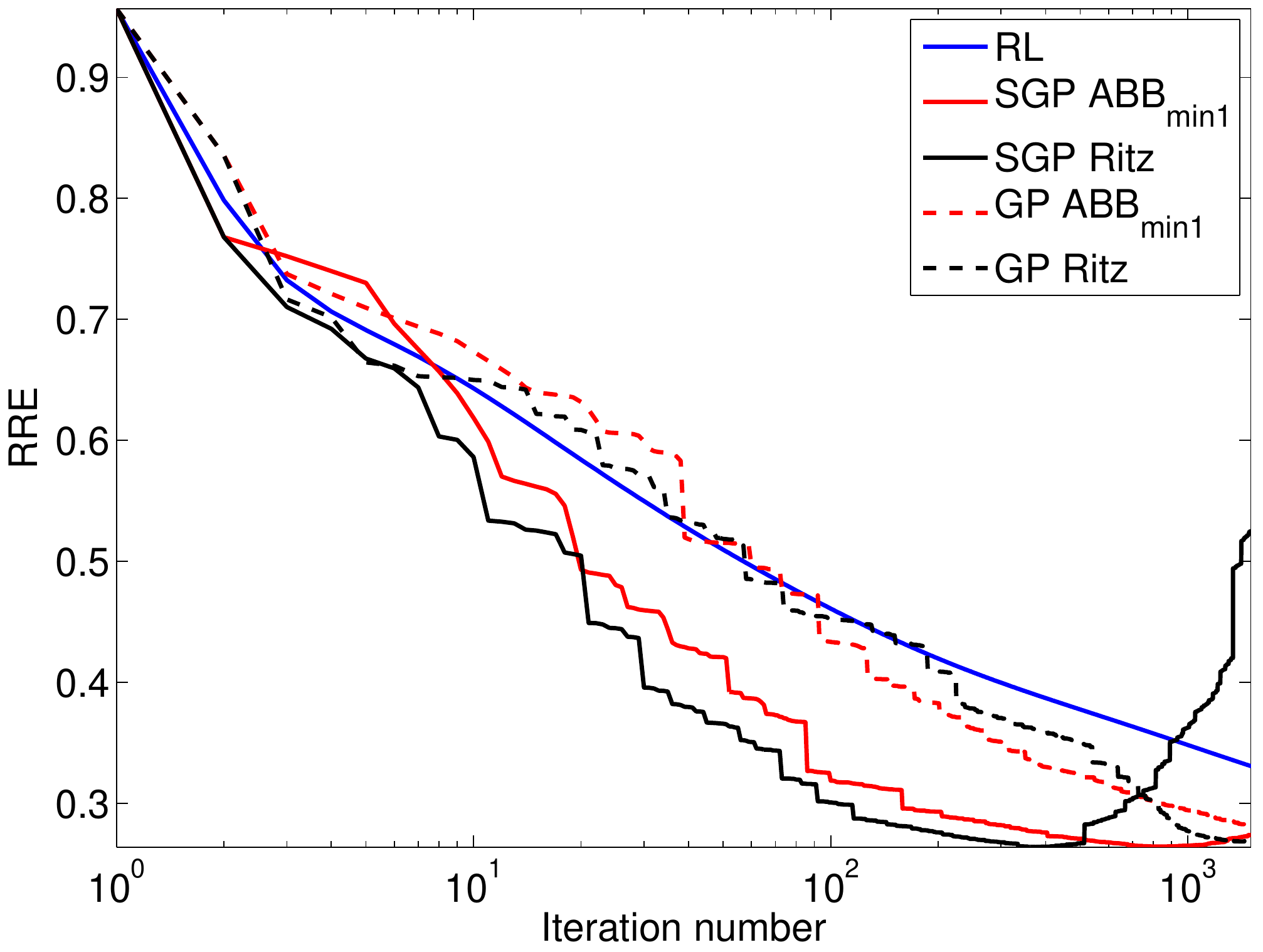}&
\includegraphics[width=.45\textwidth]{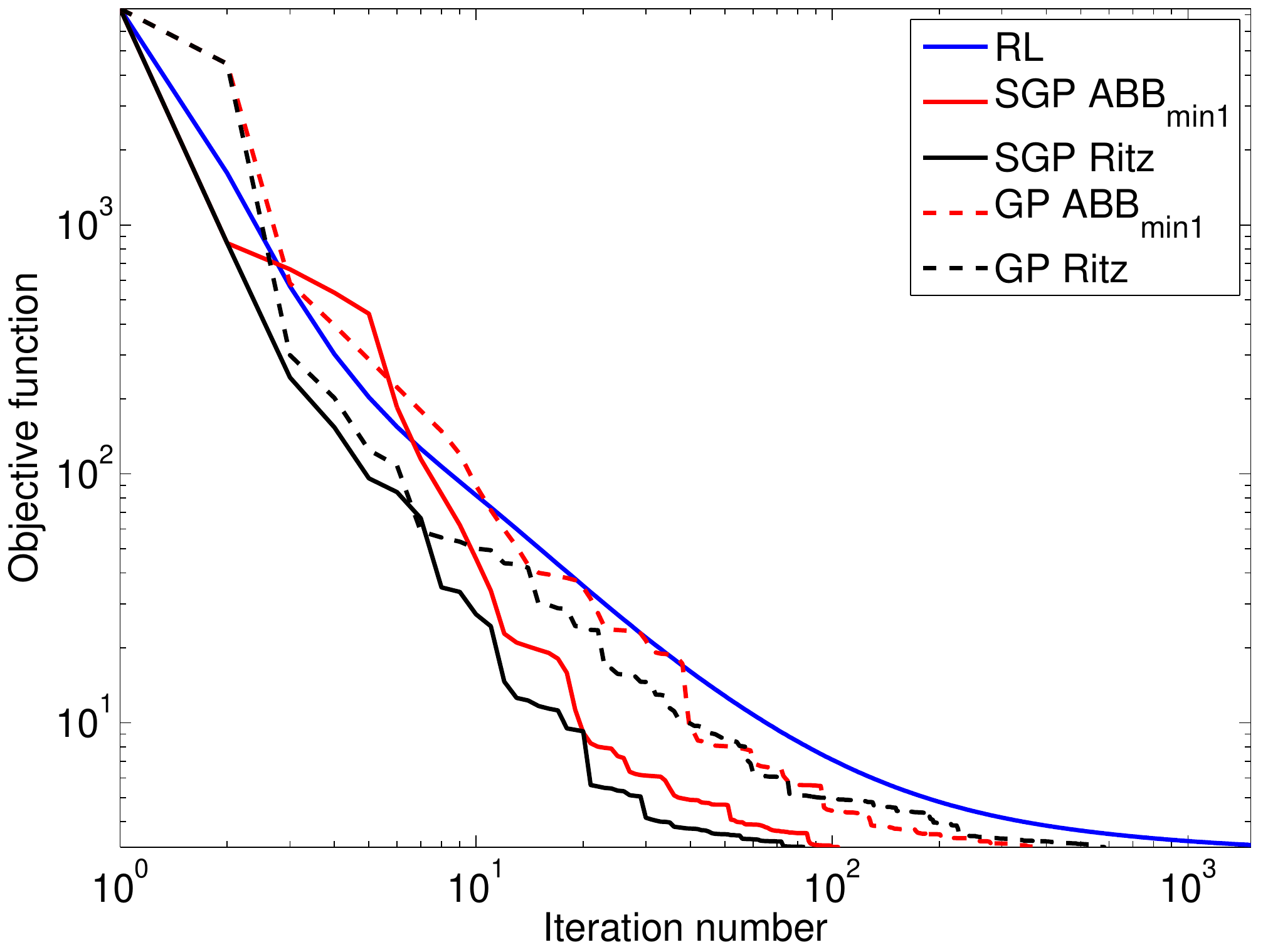}\\
\end{tabular}
\caption{Relative reconstruction error (left panel) and objective function (right panel) provided by the different methods for Image D test problem.}
\label{Err_Funct2}
\end{center}
\end{figure}

\begin{figure}[ht]
\begin{center}
\begin{tabular}{ccc}
\includegraphics[width=.3\textwidth]{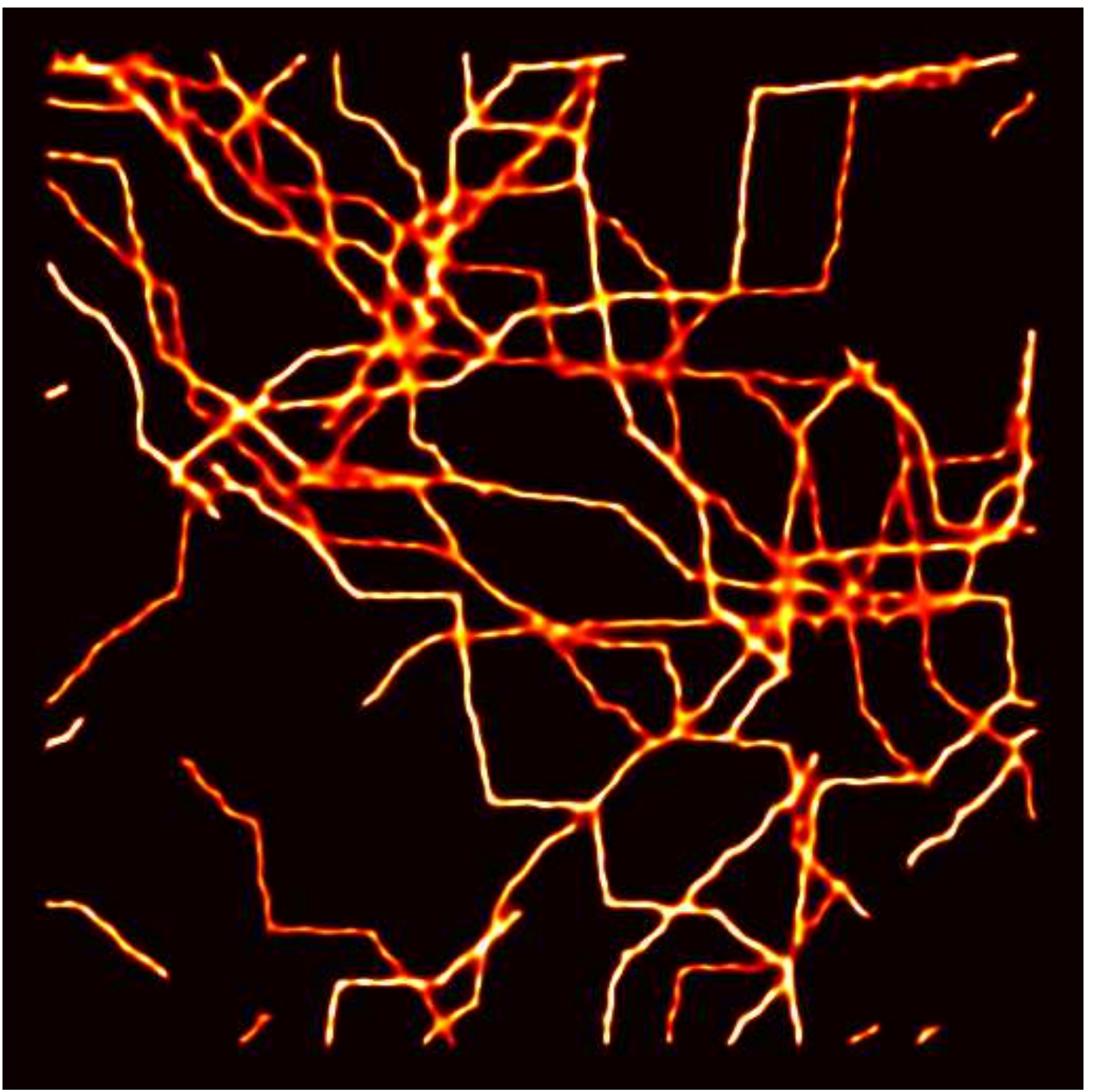}&
\includegraphics[width=.3\textwidth]{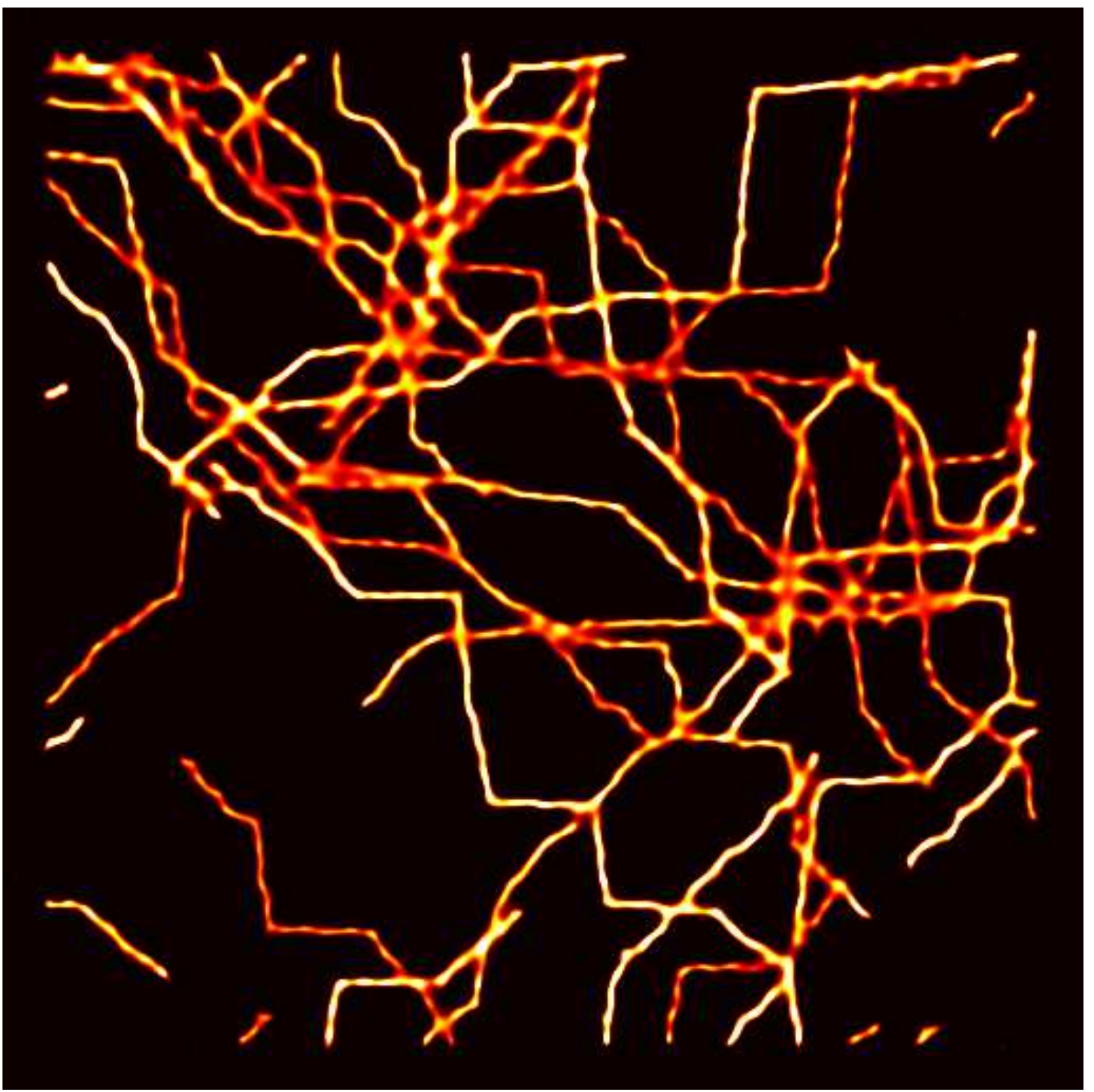}&
\includegraphics[width=.3\textwidth]{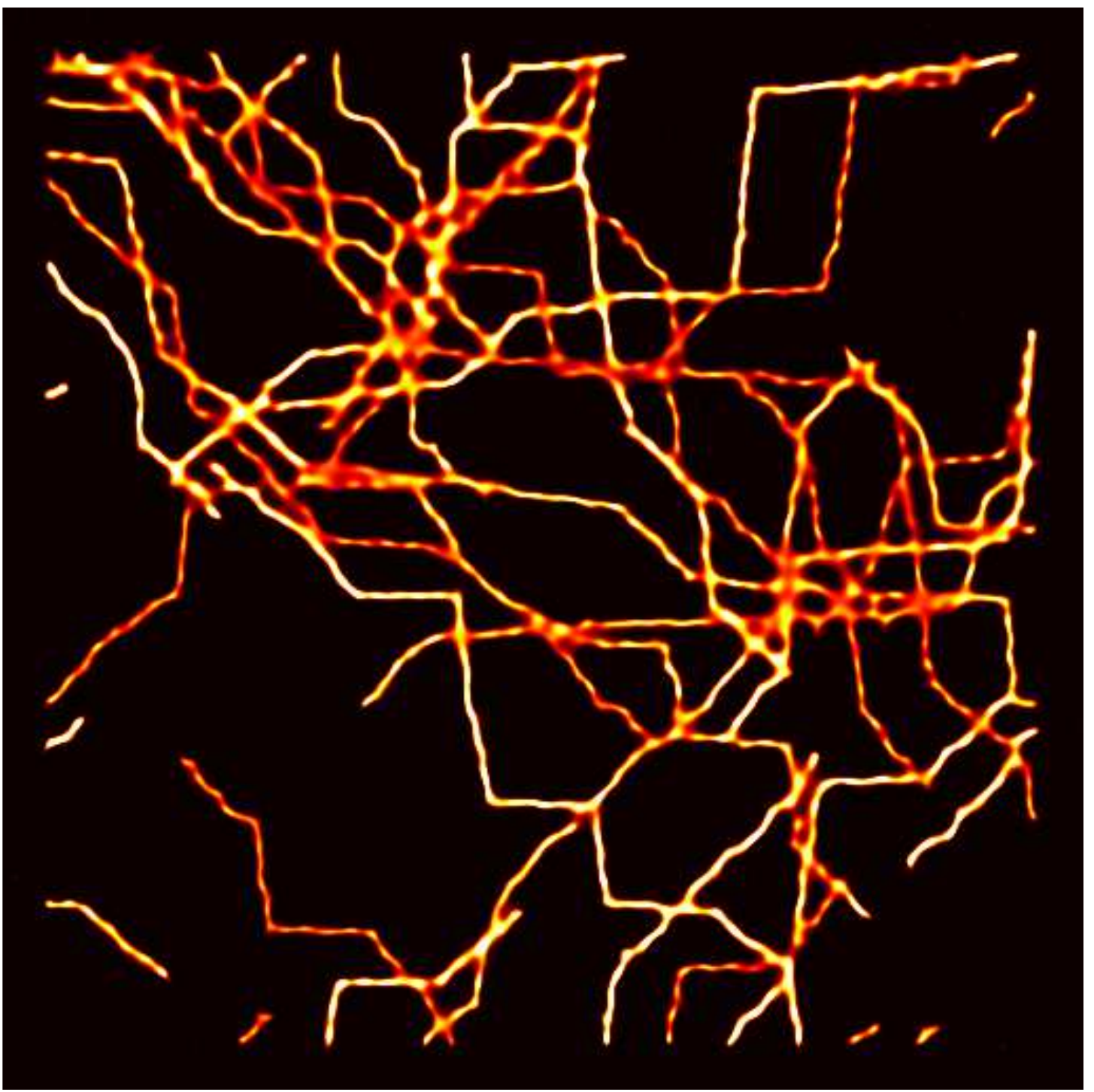}\\
RL & SGP ABB$_{\rm{min1}}$ & SGP Ritz \\
\end{tabular}
\caption{Reconstruction of the object E obtained by different SGP methods.}
\label{Rec_Poiss}
\end{center}
\end{figure}

In this section the minimization of the KL-divergence defined in \eqref{J_Poisson}, subject to non-negative constraints, on two datasets has been studied. We considered two objects of different size: the $256\times 256$ spacecraft image (used also in the Gaussian noise discussion) and a $512\times 512$ microscopy phantom representing a micro-tubule network inside the cell \cite{Porta2014}. The blurred and noisy images have been obtained by convolving the original images with the PSF described in the previous section, adding a constant background equal to 100 and 1, respectively, and by perturbing the result of the convolution with Poisson noise. Figure \ref{Immagini2} reports the images of the spacecraft and phantom datasets, indicated by D and E.\\
In our tests on Poisson data we excluded the GP Extra algorithm since a) the extrapolation step might generate a vector $\overline{\ve{x}}^{(k)}$ outside the domain of the KL divergence, and b) only an upper bound of the Lipschitz constant for $\nabla J_0^{KL}$ is available \cite{Harmany2012}. The minimum error reached by the compared methods and the corresponding number of iterations and execution time needed to recover an approximation of the true image have been reported in Table \ref{Tab2}. We also show the results obtained with the Richardson-Lucy (RL) algorithm \cite{Lucy1974,Richardson1972}, which is the strategy commonly used in the literature to treat image reconstruction problems with Poisson data and whose $(k+1)$-th iteration is defined by
\begin{equation*}
\ve{x}^{(k+1)} = {\rm{diag}}\left(\frac{\ve{x}^{(k)}}{A^T\ve{1}}\right)A^T\left(\frac{\ve{y}}{A\ve{x} + \ve{b}}\right) = \ve{x}^{(k)} - {\rm{diag}}\left(\frac{\ve{x}^{(k)}}{A^T\ve{1}}\right)\nabla J_0^{KL}(\ve{x}^{(k)}).
\end{equation*}
As shown in the previous equation, also the RL algorithm can be viewed as a scaled gradient method with constant steplength equal to 1.\\
We remark that, for all the considered methods, the main computations for each iteration are the two matrix-vector products $A\ve{x}^{(k)}$ and $A^T(\ve{y}/(A\ve{x}^{(k)} + \ve{b}))$, which require 4 FFTs if periodic boundary conditions are assumed \cite{Hansen2006}. The reconstruction error behavior and the decrease of the objective function generated by the different algorithms in solving test problem D can be appreciated in Figure \ref{Err_Funct2}, while in Figure \ref{Rec_Poiss} we report the reconstructions of object E provided by RL and the SGP methods.\\

\begin{table}[ht]
\caption{Minimum RRE achieved by each algorithm in the non regularized Poisson deblurring problems, with the corresponding number of iterations required and execution time.}\label{Tab2}
\begin{center}
\begin{tabular}{l|ccc|ccc|}
                      & \multicolumn{3}{c|}{{\bf D}} & \multicolumn{3}{c|}{{\bf E}} \\
											& It. & RRE & Time(s) & It. & RRE & Time(s)\\
\hline
GP ABB$_{\rm{min1}}$  & 3409 & 0.268 & 35.34 & 6444 & 0.436 & 360.7 \\
GP Ritz               & 1426 & 0.268 & 15.02 & 3998 & 0.436 & 224.0 \\
RL                    & 8426 & 0.291 & 67.26 & 1027 & 0.463 & 39.76 \\
SGP ABB$_{\rm{min1}}$ &  821 & 0.264 & 9.553 &  165 & 0.440 & 10.62 \\
SGP Ritz              &  375 & 0.264 & 4.382 &  104 & 0.442 & 6.569 \\
\hline
\end{tabular}
\end{center}
\end{table}

\subsection{Edge-preserving regularization}\label{subsec:KL+HS}

\begin{figure}[ht]
\begin{center}
\begin{tabular}{cc}
\includegraphics[width=.3\textwidth]{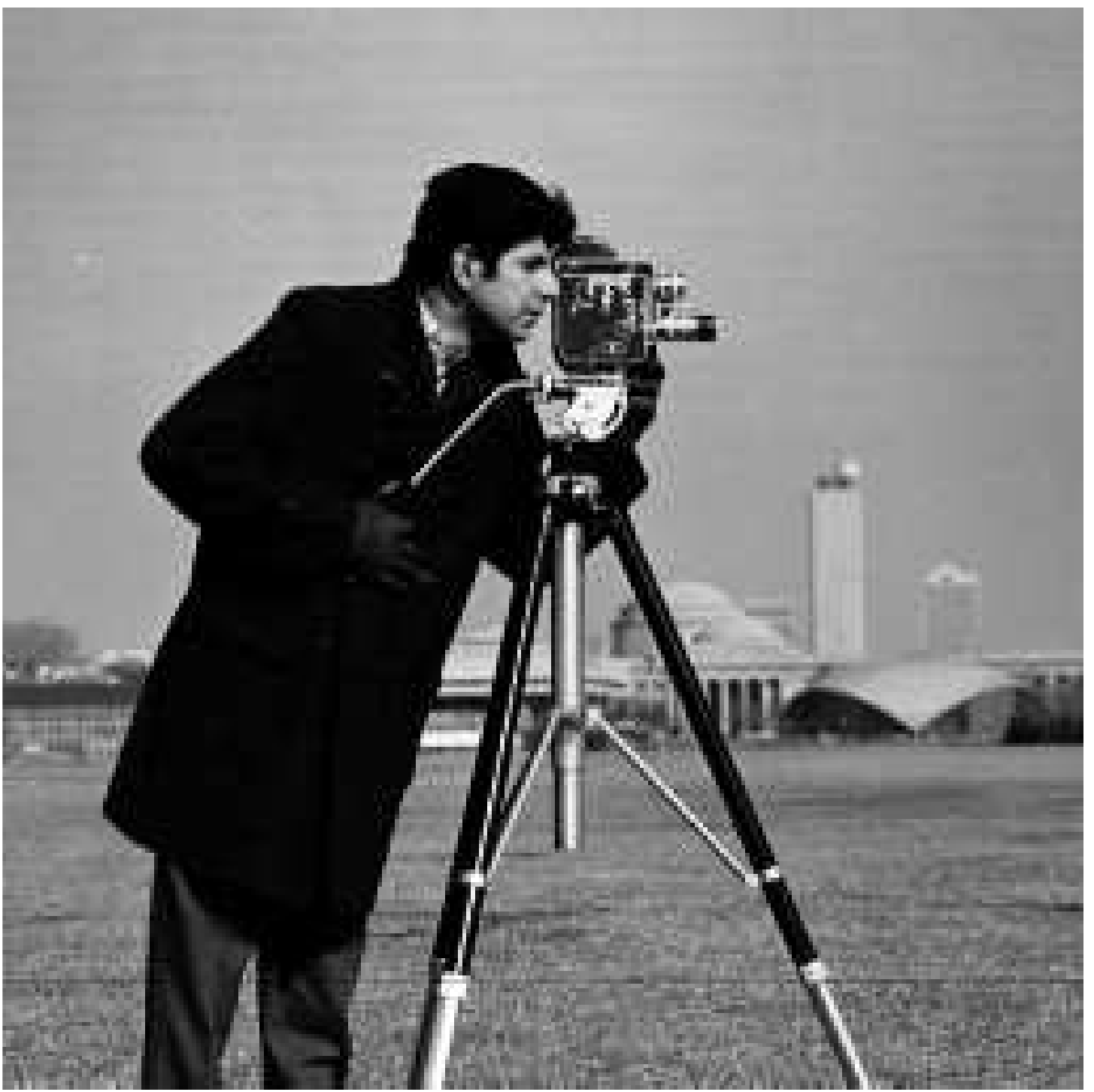}&
\includegraphics[width=.3\textwidth]{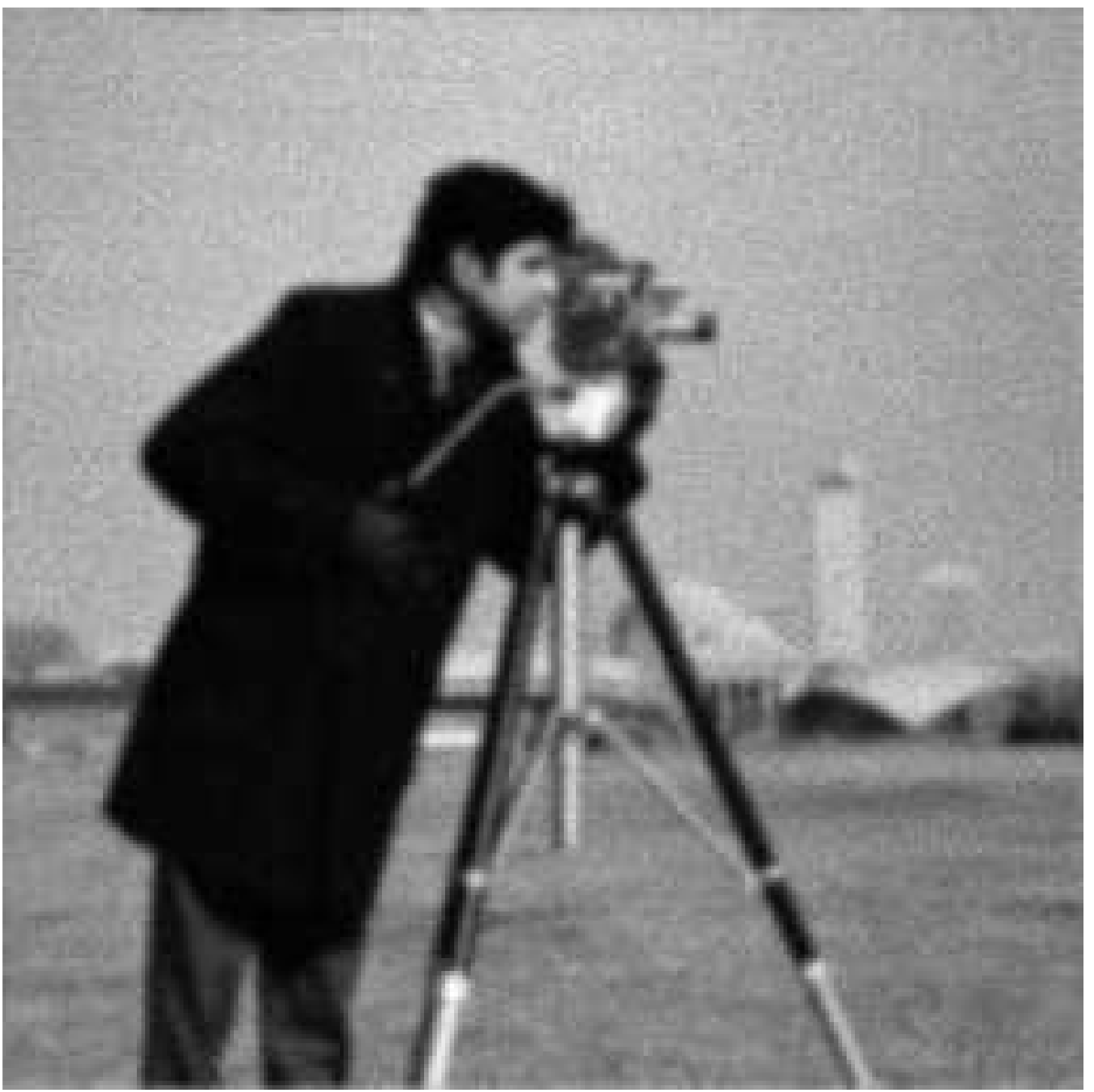}\\
\end{tabular}
\caption{The original image (left panel) and the corrupted data (right panel) for the test problem F.}
\label{Immagini3}
\end{center}
\end{figure}

\begin{figure}[ht]
\begin{center}
\includegraphics[width=.45\textwidth]{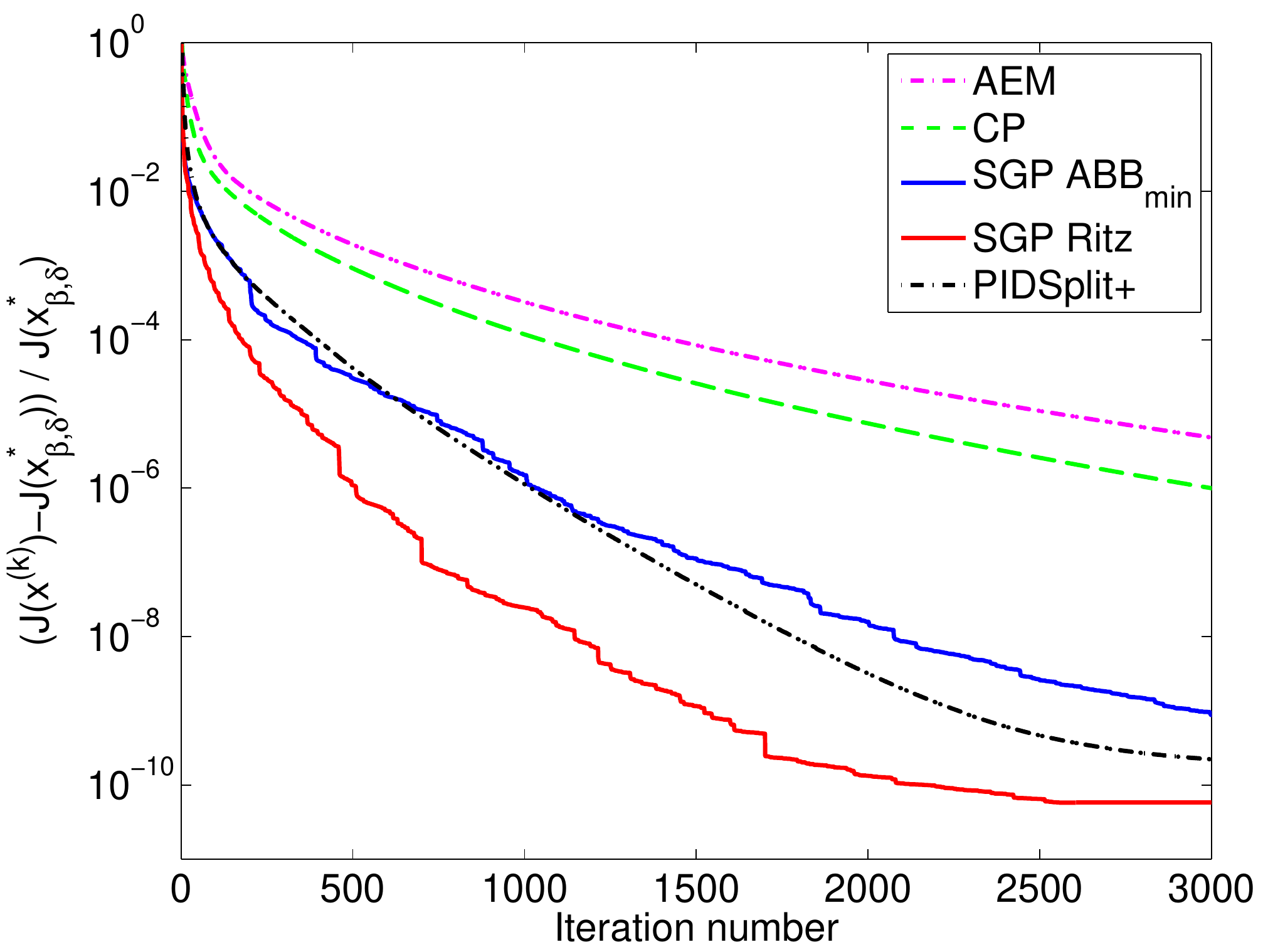}
\caption{Relative difference \eqref{RED} between the objective function $J(\ve{x}^{(k)})$ and the minimum value $J(\ve{x}^*_{\beta,\delta})$ provided by the different methods for the test problem F.}
\label{f_fmin}
\end{center}
\end{figure}

A further numerical test has been performed by solving the regularized minimization problem \eqref{Min_prob_reg} with the HS term defined in \eqref{TV_smooth} on the dataset, called F, shown in Figure \ref{Immagini3}. The original image is the $256 \times 256$ Cameraman used in several papers (see e.g. \cite{Zhu2008}). The values of the original image are in the range $[0,1000]$ and the background term has been set to zero. The corrupted data has been generated by convolving the object with a Gaussian PSF with standard deviation equal to 1.3 and adding Poisson noise. For these tests, we compared the two SGP approaches with other three recent methods, namely:
\begin{itemize}
\item[$\bullet$] the PIDSplit+ algorithm \cite{Setzer2010}, which is an alternating direction method of multipliers (ADMM) specifically tailored for the non-negative minimization of the KL functional with the addition of the total variation regularization term;
\item[$\bullet$] the alternating extragradient method (AEM) \cite{Bonettini2011b} and the Chambolle\&Pock (CP) algorithm \cite{Chambolle2011}, which are strategies for saddle point problems which apply to the minimization of a sum of convex functions reformulated in primal-dual form.
\end{itemize}
The original schemes have been suitably adapted by ourselves to account for the presence of the smoothing parameter $\delta$ in the HS regularization term.\\
We add a few remarks on the computational cost of a single iteration for the considered approaches to clarify the comparison of our results. As in the non regularized case, all the methods need the computation of two matrix-vector products involving $A$ and $A^T$. In addition, the PIDSplit+ algorithm requires the solution of a $n^2 \times n^2$ linear system, which can be computed by means of two FFTs exploiting the structure of the coefficient matrix \cite{Setzer2010}. As concerns AEM, the additional number of matrix-vector products depends on the backtracking procedure needed to set the steplength parameter. We remark that AEM is a fully automatic scheme (i.e., all its parameters are self-tuned) while CP and PIDSplit+ depends on user supplied parameters whose choice strongly influences its convergence behaviour (see e.g. \cite{Bonettini2012}).\\
After several tests we chose $\beta = 0.0045$ and $\delta=0.1$ as values of the parameters providing good reconstructions, and we performed 100000 PIDSplit+ iterations (with the value of its parameter $\gamma$ set equal to $50/\beta$ as in \cite{Setzer2010}) to get an approximate solution $\ve{x}^*_{\beta,\delta}$. Then we run AEM, CP, SGP ABB$_{\rm{min1}}$ and SGP Ritz and took note of the first iterations when the relative difference between the objective function and the minimum value
\begin{equation}\label{RED}
\frac{J(\ve{x}^{(k)}) - J(\ve{x}^*_{\beta,\delta})}{J(\ve{x}^*_{\beta,\delta})}
\end{equation}
was below certain thresholds (e.g., $10^{-4}$, $10^{-6}$ and $10^{-8}$). As concerns the CP algorithm, we tried different values of the parameters $(\tau,\sigma)$ satisfying the condition needed by the method to converge (see \cite[Theorem 1]{Chambolle2011}) and we found our best results by setting $\tau=0.001$ and $\sigma=100$. Table \ref{Tab3} shows the numbers of iterations needed together with the execution times. In all cases the corresponding reconstruction errors (i.e., the relative Euclidean errors between the $k$-th iterate and the true object) have been equal to 0.087. We remark that, when an explicit regularization term is present in the objective function, the optimization algorithms can be compared only in terms of efficiency, since the quality of the reconstruction depends only on the selected regularization term and the choice of the regularization parameter. The information on the RRE between the current iterate and the true object are provided only for sake of completeness. The plots of the distances \eqref{RED} as functions of the iterations obtained by applying AEM, CP, SGP ABB$_{\rm{min1}}$, SGP Ritz and PIDSplit+ are shown in Figure \ref{f_fmin}.\\
The presence of a HS regularization term in the objective function leads to similar conclusions as using the non regularized problems. In fact, the combination between SGP and the limited memory steplength allows again a substantial reduction of the iterations, and results to be comparable with more elaborated strategies requiring a heavier cost per iteration.

\begin{table}[ht]
\caption{Numbers of iterations and execution times required by each algorithm to bring the relative difference between the objective function and the minimum below given thresholds.}
\label{Tab3}
\begin{center}
\setlength{\tabcolsep}{4pt}
\begin{tabular}{l|cc|cc|cc|}
											& \multicolumn{6}{c|}{{\bf F}} \\
											& \multicolumn{2}{c|}{Tol = $10^{-4}$}  & \multicolumn{2}{c|}{Tol = $10^{-6}$} & \multicolumn{2}{c|}{Tol = $10^{-8}$} \\
                      & It. & Time(s) & It.  & Time(s) & It.  & Time(s) \\
\hline
AEM                   & 1428 & 83.97  & 4101 & 242.9   & 9140 & 539.1  \\
CP										& 1049 & 39.75  & 2998 & 112.4   & 6556 & 244.7  \\
SGP ABB$_{\rm{min1}}$ &  347 & 18.37  & 1032 & 61.57   & 2076 & 145.1  \\
SGP Ritz              &  179 & 11.27  &  510 & 31.23   & 1146 & 70.15  \\
PIDSplit+             &  398 & 34.02  & 1019 & 82.89   & 1783 & 143.1  \\
\hline
\end{tabular}
\end{center}
\end{table}

\section{Beyond non-negativity}\label{sec6}

\begin{figure}[ht]
\begin{center}
\begin{tabular}{cc}
\includegraphics[width=.3\textwidth]{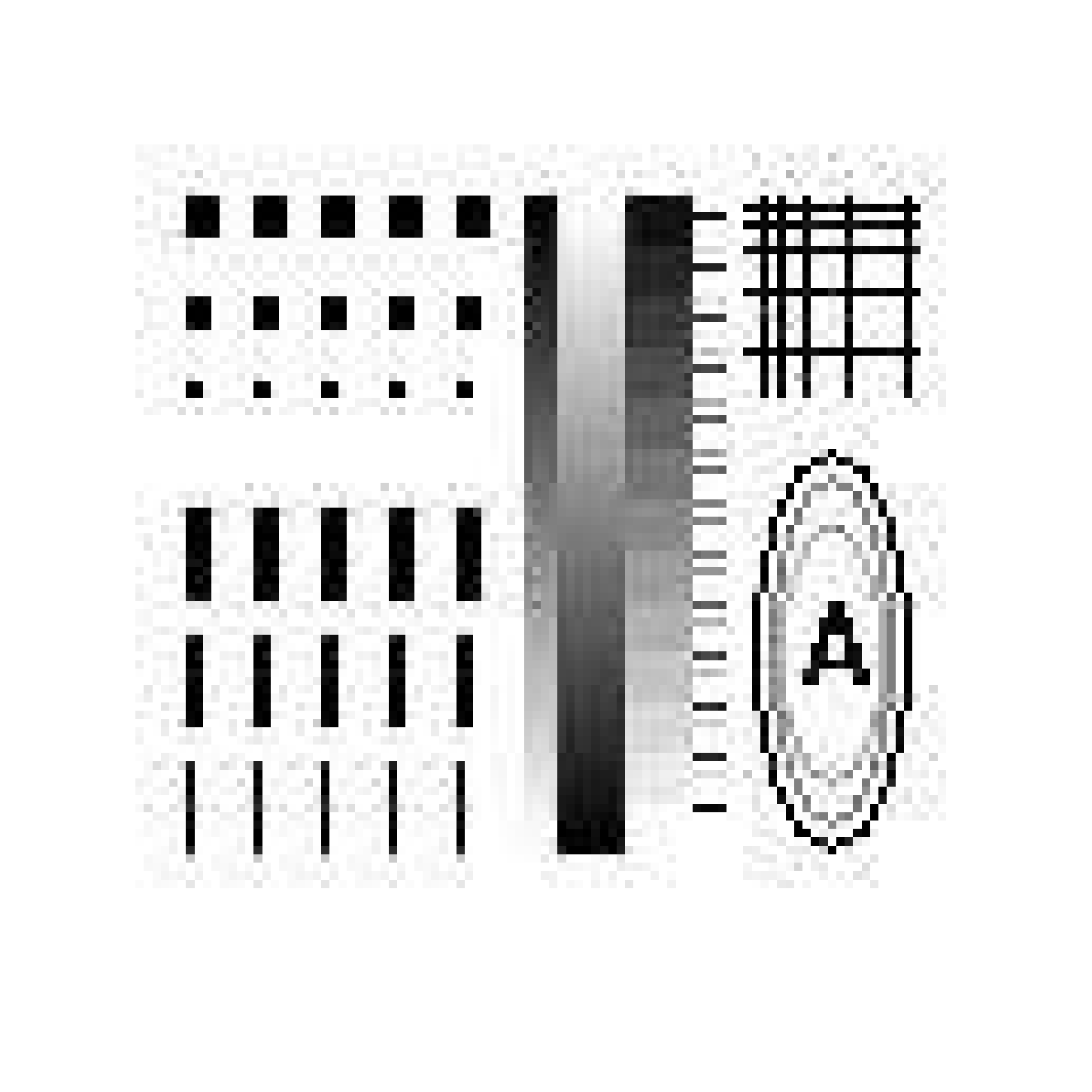}&
\includegraphics[width=.3\textwidth]{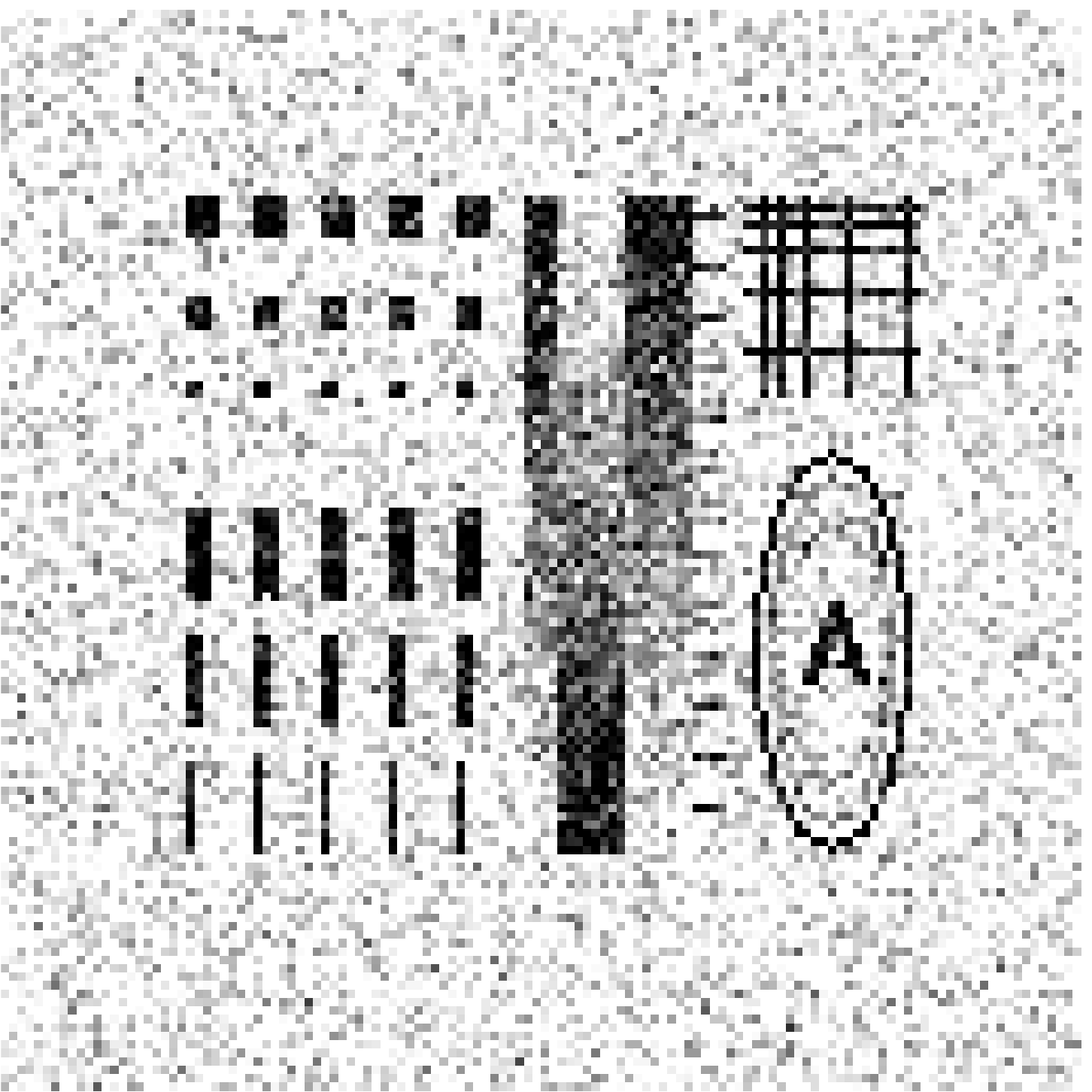}\\
\end{tabular}
\caption{The original image (left panel) and the corrupted data (right panel) for the test problem G.}
\label{Immagini4}
\end{center}
\end{figure}

\begin{figure}[ht]
\begin{center}
\includegraphics[width=.45\textwidth]{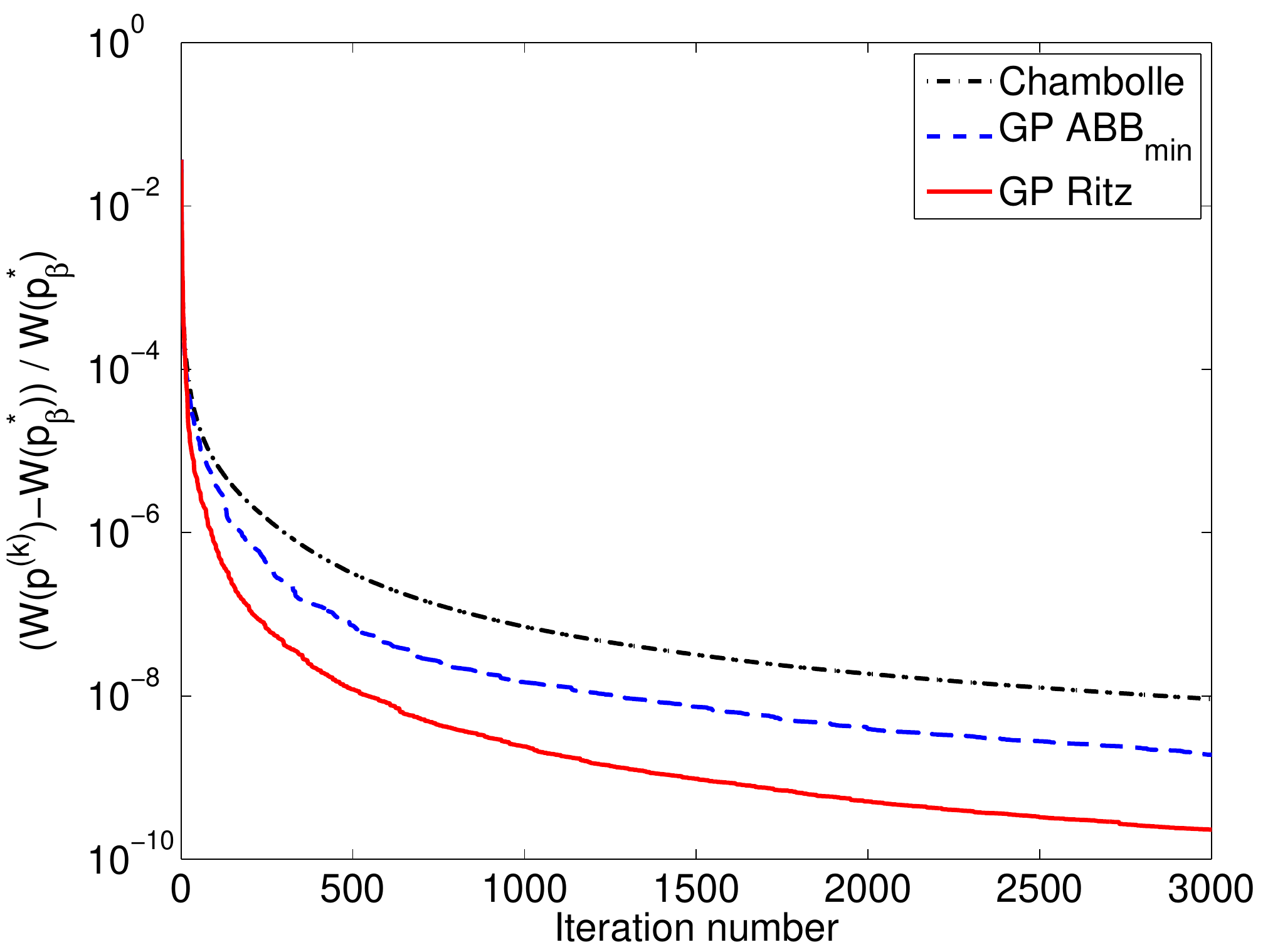}
\caption{Relative difference \eqref{RED2} between the objective function $\mathcal{W}(\ve{p}^{(k)})$ and the minimum value $\mathcal{W}(\ve{p}^*_{\beta})$ provided by the different methods for the test problem G.}
\label{f_fmin2}
\end{center}
\end{figure}

In this section we provide some hints to generalize the limited memory steplength rule described in the previous sections to different constraints. As a test problem, we consider the total variation based image denoising model proposed by Rudin, Osher and Fatemi (ROF) \cite{Rudin1992}. The discrete ROF model aims at solving the optimization problem
\begin{equation}\label{eq:PrimalROF}
\min_{\ve{x}\in\mathbb{R}^{n^2}} \frac{1}{2\beta}\|\ve{x} - \ve{y}\|_2^2 + J_R^{TV}(\ve{x})
\end{equation}
given the data $\ve{y}\in\mathbb{R}^{n^2}$ and the regularization parameter $\beta > 0$. The functional $J_R^{TV}$ denotes the discrete version of the total variation defined as in \eqref{TV_smooth} with $\delta$ equal to zero. We remark that the solution of \eqref{eq:PrimalROF} is uniquely defined as a consequence of the strict convexity of the objective function. However, the nondifferentiability of $J_R^{TV}$ prevents us from directly applying a gradient method in order to find the minimum point of \eqref{eq:PrimalROF}. According to \cite{Chambolle2004}, a strategy to overcome this difficulty consists in taking into account the dual formulation of the primal problem \eqref{eq:PrimalROF}:
\begin{equation}\label{eq:DualROF}
\min_{\ve{p} \in \mathcal{P}} \mathcal{W}(\ve{p}) \equiv \big\|\beta \, \mathrm{div}(\ve{p}) - \ve{y}\big\|^2,
\end{equation}
where $\mathcal{P} = \Big\{\ve{p} \in \mathbb{R}^{2n^2} : \sqrt{p_i^2 + p_{i+n^2}^2}\leq 1, \forall i=1,...,n^2\Big\}$ is the feasible set and the discrete divergence operator $\mathrm{div} : \mathbb{R}^{2n^2}\longrightarrow\mathbb{R}^{n^2}$ is the adjoint of the discrete gradient operator \eqref{discr_grad}.  More in details, the identity $\langle \mathcal{D} \ve{x}, \ve{p}\rangle_{\mathbb{R}^{2n^2}} = - \langle \ve{x}, \mathrm{div}(\ve{p}) \rangle_{\mathbb{R}^{n^2}}$ defines the divergence operator uniquely. The next proposition, proved in \cite{Zhu2008}, establishes how to get the primal solution starting from a dual solution.
\begin{proposizione}
 If $\{\ve{p}^{(k)}\}_{k\in\mathbb{N}}\subset \mathcal{P}$ is a sequence such that all its accumulation points are stationary points of \eqref{eq:DualROF}, then the sequence $\{\ve{x}^{(k)}\}_{k\in\mathbb{N}} = \{\ve{y} - \beta \mathrm{div}(\ve{p}^{(k)})\}_{k\in\mathbb{N}}$ converges to the unique solution of \eqref{eq:PrimalROF}.
\end{proposizione}
Therefore we address the problem of finding the stationary points of \eqref{eq:DualROF}, which are global minimum points due to the convexity of $\mathcal{W}$. The dual ROF \eqref{eq:DualROF} is a differentiable, constrained minimization problem and can be addressed with gradient projection algorithms. A very popular scheme belonging to this class has been proposed by Chambolle \cite{Chambolle2004} and is described by the iteration
\begin{equation}\label{eq:Chambolle_alg}
p_i^{(k+1)} = \frac{p_i^{(k)} - \tau \nabla \mathcal{W}(\ve{p}^{(k)})_i}{\sqrt{\big(p^{(k)}_i\big)^2 + \big(p^{(k)}_{i+n^2}\big)^2}}, \quad i=1,\ldots,n^2,
\end{equation}
where the steplength $\tau$ is constant during the iterations and fixed less than $1/4$ in order to assure the convergence.\\
Among the gradient projection methods along the feasible directions we analyzed in the previous sections and used in the numerical experiments, we consider only the nonscaled approaches GP ABB$_{\mathrm{min1}}$ and GP Ritz, since we did not find any effective strategy to design a scaling matrix for \eqref{eq:DualROF}. As concerns the latter approach, we extended the steplength selection rule described in section \ref{sec3} to the constraints set $\mathcal{P}$ by preserving only the components of the gradient $\ve{g}^{(k)}$ corresponding to the nonprojected ones of $\ve{p}^{(k+1)}$. More in details, by recalling the general iteration for a  gradient projection method \eqref{SGP_iter}, a criterion to realize this idea is achieved by looking for the indexes $j$ such that
\begin{equation}\label{eq:crit_denoising}
\ell_j^{(k)} = |\big(\ve{d}^{(k)}\big)_j + \alpha_k \big(\nabla \mathcal{W}(\ve{p}^{(k)})\big)_j| < \epsilon
\end{equation}
for $\epsilon$ sufficiently small. In this way we are able to reproduce the original limited memory scheme in the suitable subset of the nonprojected components. Accordingly, the stored vectors $\widetilde{\ve{g}}^{(k)}$ are fixed by exploiting the inequality \eqref{eq:crit_denoising}:
\begin{equation}\label{critROF}
\widetilde{g}^{(k)}_j = \begin{cases} 0 & {\rm{if}} \ \ell_j^{(k)}\geq\epsilon \\ \left[\nabla \mathcal{W}(\ve{p}^{(k)})\right]_j & {\rm{if}} \ \ell_j^{(k)}<\epsilon \end{cases},
\end{equation}
and the matrix $\widetilde{\ve{G}}$ is given by
\begin{equation*}
\widetilde{G} = \left[\widetilde{\ve{g}}^{(k-m)}, \ldots, \widetilde{\ve{g}}^{(k-1)}\right].
\end{equation*}
If applied to non-negative constraints, the criterion in \eqref{critROF} to select the components of the gradients to be preserved is equivalent to the one in \eqref{critnonneg}.\\
The GP ABB$_{\mathrm{min1}}$, GP Ritz and Chambolle methods have been tested in the ROF model applied to the $128 \times 128$ Shape image available in Wright's TV-Regularized Image Denoising Software \cite{Zhu2008} and corrupted with additive white Gaussian noise with variance 1 (this dataset will be denoted by G and is shown in Figure \ref{Immagini4}). The regularization parameter $\beta$ has been selected equal to $20$.\\
As done in section \ref{subsec:KL+HS}, we compare the performances of the different methods by checking how well they approach an approximate solution $\ve{p}^*_{\beta}$ of \eqref{eq:DualROF}, computed by running GP Ritz for 100000 iterations. In particular, Table \ref{Tab4} reports the number of iterations and the execution time needed by the considered algorithms to bring the relative difference between the objective function and the minimum value
\begin{equation}\label{RED2}
\frac{\mathcal{W}(\ve{p}^{(k)}) - \mathcal{W}(\ve{p}^*_{\beta})}{\mathcal{W}(\ve{p}^*_{\beta})}
\end{equation}
less than certain thresholds (e.g., $10^{-4}$, $10^{-6}$ and $10^{-8}$).
In all cases, the corresponding reconstruction errors on the primal solution (i.e., the relative Euclidean errors between the $k$-th approximation $\ve{x}^{(k)}$ and the true object) have been equal to 0.128. Figure \ref{f_fmin2} shows the distances defined in \eqref{RED2} and provided by the Chambolle, GP ABB$_{\rm{min1}}$ and GP Ritz methods versus the iteration number.\\

\begin{table}[ht]
\caption{Numbers of iterations and execution times required by each algorithm to bring the relative difference between the objective function and the minimum below given thresholds.}
\label{Tab4}
\begin{center}
\setlength{\tabcolsep}{4pt}
\begin{tabular}{l|cc|cc|cc|}
											& \multicolumn{6}{c|}{{\bf G}} \\
											& \multicolumn{2}{c|}{Tol = $10^{-4}$}  & \multicolumn{2}{c|}{Tol = $10^{-6}$} & \multicolumn{2}{c|}{Tol = $10^{-8}$} \\
                      & It. & Time(s) & It.  & Time(s) & It.  & Time(s) \\
\hline
CP										&  19 & 0.343  & 301 & 4.804   & 2854 & 44.36  \\
GP ABB$_{\rm{min1}}$ &  14 & 0.499  & 177 & 5.865   & 1263 & 40.70  \\
GP Ritz              &  15 & 0.702  &  89 & 3.369   &  543 & 19.56  \\
\hline
\end{tabular}
\end{center}
\end{table}

The results presented both in Table \ref{Tab4} and in Figure \ref{f_fmin2} confirm the goodness of the suggested limited memory steplength selection scheme for a gradient projection method with respect to standard approaches also for a constrained optimization problem where the feasible set is different from the simple non-negative orthant.

\section{Conclusions}\label{sec7}

In this paper we considered a first-order method for the minimization of non-negatively constrained optimization problems arising in the image reconstruction field, and we introduced a new strategy for the steplength selection which generalizes a rule recently proposed in the unconstrained optimization framework. The steplength value is based on the storage of a limited number of consecutive objective function gradients and we showed how it can be extended to account for the presence of both a scaling matrix multiplying the gradient of the objective function and a non-negative constraint on the pixels of the unknown image. We first tested our rule in the minimization of a quadratic function with different features, and we showed that the limited memory steplength is extremely competitive with respect to state-of-the-art BB-like choices. Similar conclusions can be drawn by the numerical experiments we carried out on image reconstruction problems where the measured images are affected by either Gaussian or Poisson noise. A final test on the ROF model showed the potentiality of the proposed rule also in optimization problems with different constraints.\\
Thanks to the significant reduction of the iterations achievable by the proposed steplength, in our future work we will consider the application of our new scheme to real-world imaging problems, as the reconstruction of X-ray images of solar flares starting from the emitted radiation \cite{Bonettini2010a,Bonettini2014} and the deblurring of conventional stimulated emission depletion (STED) microscopy images of sub-cellular structures in fixed cells \cite{Zanella2013b}. Moreover, the proposed rule will be tested also within a SGP method where the sequence of scaling matrices converges to the identity, since in this case strong convergence results have been recently proved under mild convexity assumptions \cite{Bonettini2015}.

\section*{Acknowledgments}
This work has been partially supported by the Italian Spinner 2013 PhD Project ``High-complexity inverse problems in biomedical applications and social systems'' and by MIUR (Italian Ministry for University and Research), under the projects FIRB - Futuro in Ricerca 2012, contract RBFR12M3AC, and PRIN 2012, contract 2012MTE38N. The Italian GNCS - INdAM (Gruppo Nazionale per il Calcolo Scientifico - Istituto Nazionale di Alta Matematica) is also acknowledged.

\end{document}